\documentclass[11pt]{article}
\usepackage[english]{babel}
\usepackage{indentfirst}
\usepackage{latexsym}
\usepackage[colorlinks=true, bookmarksnumbered=true, bookmarksopen=true,
bookmarksopenlevel=3, pdfstartview=FitH, linkcolor=cyan, pdfmenubar=true,
pdftoolbar=true, bookmarks=true,citecolor=cyan, urlcolor=magenta,
filecolor=cyan,menucolor=black,plainpages=false,pdfpagelabels]{hyperref}
\usepackage[lmargin=2cm,tmargin=2cm,bmargin=2cm,rmargin=2cm]{geometry}
\usepackage{amsbsy, amsmath,amsfonts,amssymb,amsthm}
\usepackage{times}
\usepackage{graphicx}
\usepackage{amsmath}
\usepackage{amsfonts}
\usepackage{amssymb}
\usepackage[english]{babel}
\usepackage[latin1]{inputenc}
\usepackage[T1]{fontenc}
\usepackage{amsmath,amssymb,graphics,amsthm}
\usepackage[usenames]{color}
\usepackage{mathtools}
\usepackage{comment}

\newtheorem{Assumption}{Assumption}
\newtheorem{lema}{Lemma}[section]
\newtheorem{obs}{Remark}[section]
\newtheorem{teo}{Theorem}[section]
\newtheorem{cor}{Corollary}[section]
\newtheorem{defin}{Definition}[section]

\numberwithin{equation}{section}

\begin{document}

\title{Global well-posedness and asymptotic behavior in Besov-Morrey spaces
for chemotaxis-Navier-Stokes fluids}
\author{ {Lucas C. F. Ferreira$^{1}$}{\thanks{{Corresponding author. }%
\newline
{E-mail adresses: lcff@ime.unicamp.br (L.C.F. Ferreira),
alvesmonisse@gmail.com (M. Postigo).}}} \ \ \ and \ \ {Monisse Postigo$^{2}$}
\\
{\small $^{1,2}$Universidade Estadual de Campinas, IMECC-Departamento de
Matemática,} \\
{\small {Rua Sérgio Buarque de Holanda, 651, CEP 13083-859, Campinas-SP,
Brazil.}}}
\date{}
\maketitle

\begin{abstract}
In this work we consider the Keller-Segel system coupled with Navier-Stokes
equations in $\mathbb{R}^{N}$ for $N \geq 2$. We prove the global
well-posedness with small initial data in Besov-Morrey spaces. Our initial
data class extends previous ones found in the literature such as that
obtained by Kozono-Miura-Sugiyama (J. Funct. Anal. 2016). It allows to
consider initial cell density and fluid velocity concentrated on smooth
curves or at points depending on the spatial dimension. Self-similar
solutions are obtained depending on the homogeneity of the initial data and
considering the case of chemical attractant without degradation rate.
Moreover, we analyze the asymptotic stability of solutions at infinity and
obtain a class of asymptotically self-similar ones. \medskip


\end{abstract}


\renewcommand{\abstractname}{Abstract}

\section{Introduction}

We consider the double chemotaxis-Navier-Stokes equations in the whole space
$\mathbb{R}^{N}$
\begin{equation}
\left\{
\begin{array}{ll}
\displaystyle\partial _{t}n+u\cdot \nabla n=\Delta n-\nabla \cdot (n\nabla
c)-\nabla \cdot (n\nabla v), & \mbox{ in
}\mathbb{R}^{N}\times (0,\infty ), \\
\displaystyle\partial _{t}c+u\cdot \nabla c=\Delta c-nc, & \mbox{ in }%
\mathbb{R}^{N}\times (0,\infty ), \\
\displaystyle\partial _{t}v+u\cdot \nabla v=\Delta v-\gamma v+n, &
\mbox{ in
}\mathbb{R}^{N}\times (0,\infty ), \\
\displaystyle\partial _{t}u+(u\cdot \nabla )u=\Delta u-\nabla \pi -nf, &
\mbox{
in }\mathbb{R}^{N}\times (0,\infty ), \\
\nabla \cdot u=0, & \mbox{ in }\mathbb{R}^{N}\times (0,\infty ), \\
n(x,0)=n_{0}(x),\,\,\,c(x,0)=c_{0}(x),\,\,\,v(x,0)=v_{0}(x),\,\,%
\,u(x,0)=u_{0}(x), & \mbox{ in }\mathbb{R}^{N}, \\
&
\end{array}%
\right.  \label{NSC}
\end{equation}%
where $N\geq 2$ and $\gamma \geq 0.$ The unknown $n(x,t),c(x,t),$ $%
v(x,t),u(x,t)$ and $\pi (x,t)$ stands for cell density, oxygen
concentration, chemical-attractant concentration, fluid velocity field, and
pressure of the fluid, respectively. The time-independent field $f$ denotes
a force field acting on the motion of the fluid.

The system (\ref{NSC}) was introduced by Tuval \textit{et al}. in \cite%
{Tuval} and corresponds to a double chemotaxis model that describes the
movement of swimming bacteria living in an incompressible viscous fluid,
which swim toward a higher concentration of oxygen and chemical attractant.
The fluid movement is modeled by the Navier-Stokes equations under the
influence of a force $-nf$ that can be produced by different mechanisms,
e.g., force due to the aggregation of bacteria onto the fluid generating a
buoyancy-like force. In turn, the chemical attractant $v$ is produced by the
bacteria themselves that degrades at a constant rate $\gamma \geq 0$.

In the case with no chemical-attractant degradation rate (i.e., $\gamma =0$)$%
,$ system (\ref{NSC}) has the scaling (see Section \ref{sec3})
\begin{equation}
(n,c,v,u)\rightarrow (\lambda ^{2}\,n(\lambda x,\lambda ^{2}t),c(\lambda
x,\lambda ^{2}t),v(\lambda x,\lambda ^{2}t),\lambda \,u(\lambda x,\lambda
^{2}t)),  \label{scal1}
\end{equation}%
which, by taking $t=0,$ induces the initial data scaling%
\begin{equation}
(n_{0},c_{0},v_{0},u_{0})\rightarrow (\lambda ^{2}\,n(\lambda x),c(\lambda
x),v(\lambda x),\lambda \,u(\lambda x)).  \label{scal2}
\end{equation}%
For mathematical analysis purposes, the above scaling relations also work
well for the case $\gamma \neq 0$ and functional spaces invariant by them
are called critical ones for (\ref{NSC}). Throughout the paper, spaces of
scalar and vector functions are denoted in the same way, e.g., we write $%
u_{0}\in L^{p}(\mathbb{R}^{N})$ in place of $u_{0}\in (L^{p}(\mathbb{R}%
^{N}))^{N}$.

In what follows we give a review of some results about (\ref{NSC}) and
related systems. Firstly we recall the classical Keller-Segel system
(without fluid coupling)
\begin{equation}
\left\{
\begin{array}{l}
\partial _{t}n=\Delta n-\nabla \cdot (n\nabla v), \\
\varepsilon \partial _{t}v=\Delta v-\gamma v+\beta n,%
\end{array}%
\right.  \label{KS-aux-1}
\end{equation}%
which has been studied by several authors. It models aggregation of
biological species (e.g. amoebae and bacteria) moving towards high
concentration of a chemical secreted by themselves or of food molecules
(e.g. glucose). We have the two basic cases $\varepsilon =0$
(parabolic-elliptic) and $\varepsilon =1$ (parabolic-parabolic). For the
former with $\gamma =0$, it is well known that there is a threshold value $%
8\pi /\beta $ for the initial mass that decides between global existence and
the finite time blow-up. For further details, see e.g.
Blanchet-Dolbeault-Perthame \cite{Blanchet} and Dolbeault-Perthame \cite%
{Dolbeault}. Considering the parabolic-parabolic case and a 2D smooth
bounded domain $\Omega $ with Neumann conditions, Nagai-Senba-Yoshida \cite%
{Nagai} proved global-in-time existence of solutions for $\varepsilon
=1,\gamma \geq 0$, $\beta >0$ and nonnegative initial data $n_{0},v_{0}\in
H^{1+\delta }(\Omega )$ with mass $M=\int_{\Omega }n_{0}<4\pi /\beta $.%
\textbf{\ }Horstman \cite{Horstman3} considered (\ref{KS-aux-1}) with $\beta
(n-1)$ in place of $\beta n$ and showed the existence of blow-up solutions
with $M=\left\vert \Omega \right\vert >4\pi /\beta $ and $\beta \left\vert
\Omega \right\vert \neq 4m\pi ,$ $m\in \mathbb{N},$ and without assuming any
symmetry properties. In the radial setting, the mass $8\pi /\beta $ is a
threshold value for the existence or blow-up of solutions (see \cite%
{Herrero-Velazquez-1, Horstman2}). In the whole plane $\mathbb{R}^{2},$
Calvez-Corrias \cite{Calvez-Corrias-1}\ obtained global solutions for
subcritical masses $M<8\pi /\beta ,$\ $\varepsilon =1,$ $\gamma \geq 0$ and $%
\beta >0,$\ as well as blow-up of solutions for $\varepsilon =0$ and $M>8\pi
/\beta .$ They also conjectured blow-up of solutions for $\varepsilon =1$
and $M>8\pi /\beta $ which, to the best of our knowledge, is still open. For
$N\geq 3$, Corrias-Perthame \cite{Corrias1} proved the existence of weak
solution for (\ref{KS-aux-1}) with $\varepsilon =\alpha =\beta =1$ and small
data $n_{0}\in L^{a}(\mathbb{R}^{N})$ and $\nabla v_{0}\in L^{a}(\mathbb{R}%
^{N})$ with $N/2<a\leq N$. After, Kozono-Sugiyama \cite{KozonoSugiyama-2}
considered (\ref{KS-aux-1}) with $\varepsilon =\beta =1$ and showed that if $%
\max \{1,N/4\}<a\leq N/2$ and $(n_{0},v_{0})\in H^{\frac{N}{a}-2,a}(\mathbb{R%
}^{N})\times H^{\frac{N}{a},a}(\mathbb{R}^{N})$ is small enough, then there
exists a unique global mild solution. Kozono-Sugiyama \cite{KozonoSugiyama-1}
extended the results of \cite{KozonoSugiyama-2} to $L^{(N/2,\infty )}\times
BMO$ and $N\geq 2,$ where $BMO$ stands for the space of bounded mean
oscillation functions and $L^{(p,\infty )}$ is the weak-$L^{p}$. Results on
global mild solutions for (\ref{KS-aux-1}) with small initial data belonging
to larger critical spaces can be found in the literature, namely
Besov-spaces $\dot{B}_{p,\infty }^{-(2-\frac{N}{p})}\times \dot{B}_{p,\infty
}^{\frac{N}{p}}$ \cite{Zhai1} (with $\gamma =0$ and $\varepsilon =\beta =1$%
), Morrey spaces \cite{Biler1, Wakab} and Besov-Morrey spaces $\mathcal{N}_{{%
q},{q_{1}},{\infty }}^{-(2-\frac{N}{q})}\times \dot{B}_{\infty ,\infty }^{0}$
\cite{Fer-Pre}.

In the context of chemotaxis-fluids, Duan-Lorz-Markowich \cite{Duan}
considered the model in 3D
\begin{equation}
\left\{
\begin{array}{l}
\displaystyle\partial _{t}n+u\cdot \nabla n=\delta \Delta n-\nabla \cdot
(\chi (c)n\nabla c), \\
\displaystyle\partial _{t}c+u\cdot \nabla c=\mu \Delta c-\kappa (c)n, \\
\displaystyle\partial _{t}u+(u\cdot \nabla )u=\upsilon \Delta u-\nabla \pi
-nf, \\
\nabla \cdot u=0, \\
\end{array}%
\right.  \label{NSC10}
\end{equation}%
and showed the global existence of classical solutions provided that the
initial data $(n_{0},c_{0},u_{0})$ is a small smooth perturbation of the
constant state $(n_{\infty },0,0)$ with $n_{\infty }\geq 0$. Lorz \cite{Lorz}
considered (\ref{KS-aux-1}) with $\varepsilon =0$ coupled to Stokes
equations in 2D, without oxygen concentration (i.e., $c=0$), and showed
global-in-time existence of solutions for small initial data $u_{0}\in L^{%
\frac{3}{2}}(\mathbb{R}^{2})$ and $n_{0}\in L^{1}(\mathbb{R}^{2}).$ In turn,
considering smooth bounded domains, Winkler \cite{Winkler} analyzed (\ref%
{NSC}) without chemical attractant (i.e., $v=0$) and showed the existence
and uniqueness of global classical solutions in 2D. In 3D, he considered (%
\ref{NSC}) with the evolution Stokes equation (in place of the Navier-Stokes
one) and obtained the existence of a global weak solution for $%
(n_{0},c_{0},u_{0})\in C^{0}(\overline{\Omega })\times W^{1,q}(\Omega
)\times D((-\mathbb{P}\Delta )^{s})$ where $q>n$, $s\in (3/4,1)$ and $%
\mathbb{P}$ is the Leray-Helmholz projection. Zhang \cite{Zhang} obtained
the local well-posedness for (\ref{NSC}) with $v=0$ and initial data in the
nonhomogeneous Besov spaces
\begin{equation*}
(n_{0},c_{0},u_{0})\in B_{p,r}^{s}(\mathbb{R}^{N})\times B_{p,r}^{s+1}(%
\mathbb{R}^{N})\times B_{p,r}^{s+1}(\mathbb{R}^{N})
\end{equation*}%
with $1<p<\infty $, $1\leq r<\infty $ and $s>\frac{N}{p}+1$ for $N=2,3$.
Considering (\ref{NSC}) with $v=0$, Choe-Lkhagvasuren \cite{ChoeLkha} showed
the existence of global mild solution for small initial data $%
(n_{0},c_{0},u_{0})$ in the critical Besov spaces
\begin{equation*}
\dot{B}_{r,1}^{-2+\frac{3}{r}}(\mathbb{R}^{3})\times \dot{B}_{r,1}^{\frac{3}{%
r}}(\mathbb{R}^{3})\times \dot{B}_{r,1}^{-1+\frac{3}{r}}(\mathbb{R}^{3})%
\text{ for }r\in \lbrack 1,3).
\end{equation*}%
After, Zhao-Zhou \cite{ZhaoZhou} extended the result obtained in \cite%
{ChoeLkha} to $r\in \lbrack 1,6)$. Finally, for results about
chemotaxis-fluid models with logistic terms, we refer the reader to \cite%
{Braukho},\cite{Fer-Vil},\cite{Lankeit} and their references.

In comparison with the aforementioned models, system (\ref{NSC}) consists in
a double chemotaxis-fluid model that includes the effect of both oxygen
concentration and chemical attractant. In \cite{Kozono1},
Kozono-Miura-Sugiyama obtained existence of global mild solutions for (\ref%
{NSC}) by considering $N\geq 3$, small initial data $n_{0}\in L_{w}^{\frac{N%
}{2}}$, $c_{0}\in L^{\infty }$ with $\nabla c_{0}\in L_{w}^{N}$, $v_{0}\in
\mathcal{S^{\prime }}/\mathcal{P}$ with $\nabla v_{0}\in L_{w}^{N}$ and $%
u_{0}\in L_{w}^{N},$ and small force $f\in L_{w}^{N},$ where $\mathcal{P}$
denotes the set of polynomials with $N$ variables. In the case $N=2,$ the
condition $n_{0}\in L_{w}^{1}$ is replaced by $n_{0}\in L^{1}$.

The main aim of this paper is to prove the global well-posedness for (\ref%
{NSC}) with small initial data in a larger critical framework based on
Besov-Morrey spaces $\mathcal{N}_{{q},{q_{1}},{\infty }}^{s}\left( \mathbb{R}%
^{N}\right) $ (see (\ref{BesovMorrey})). More precisely, we consider the
following critical initial data class
\begin{eqnarray}
&&n_{0}\in \mathcal{N}_{{q},{q_{1}},{\infty }}^{\frac{N}{q}-2}\left( \mathbb{%
R}^{N}\right) ,\,\,\,\,c_{0}\in L^{\infty }\left( \mathbb{R}^{N}\right)
\text{ with }\nabla c_{0}\in \mathcal{N}_{{r},{r_{1}},{\infty }}^{\frac{N}{r}%
-1}\left( \mathbb{R}^{N}\right) ,  \label{initial-data-aux-1} \\
&&v_{0}\in \mathcal{S}^{\prime }/\mathcal{P}\text{ with }\nabla v_{0}\in
\mathcal{N}_{{r},{r_{1}},{\infty }}^{\frac{N}{r}-1}\left( \mathbb{R}%
^{N}\right) ,\text{ and }u_{0}\in \mathcal{N}_{{p},{p_{1}},{\infty }}^{\frac{%
N}{p}-1}\left( \mathbb{R}^{N}\right) ,  \notag
\end{eqnarray}%
where the exponents $p$, $p_{1}$, $q$, $q_{1}$, $r$, $r_{1}$ and $N_{1}$
satisfy suitable conditions (for more details, see Assumption \ref%
{Assumption} in Section \ref{sec3}). For the above exponents, we have the
strict continuous inclusions
\begin{eqnarray}
L^{1} &\hookrightarrow &\mathcal{M}\hookrightarrow \mathcal{N}_{{q},{q_{1}},{%
\infty }}^{\frac{2}{q}-2}\,(N=2),\,\,\,\,L_{w}^{\frac{N}{2}}\hookrightarrow
\mathcal{M}_{\frac{N}{2}\frac{q_{1}}{q}}^{\frac{N}{2}}\hookrightarrow
\mathcal{N}_{{q},{q_{1}},{\infty }}^{\frac{N}{q}-2}\text{ }(q_{1}<q),\,\,\,\,
\label{imer-aux-1} \\
L_{w}^{N} &\hookrightarrow &\mathcal{M}_{N\frac{r_{1}}{r}}^{N}%
\hookrightarrow \mathcal{N}_{{r},{r_{1}},{\infty }}^{\frac{N}{r}-1}\text{ }%
(r_{1}<r),\,\text{\ and }\,\,L_{w}^{N}\hookrightarrow \mathcal{N}_{{p},{p_{1}%
},{\infty }}^{\frac{N}{p}-1}\text{ }(p_{1}<p),  \notag
\end{eqnarray}%
where $\mathcal{M}_{p_{1}}^{p}$ denotes Morrey spaces and $\mathcal{M}%
_{1}^{1}=\mathcal{M}$ the space of finite signed Radon measures (see %
\eqref{morrey}). In particular, depending on the spatial dimension $N$ and
taking suitable values for the indexes in (\ref{imer-aux-1}), we can
consider the initial cell density $n_{0}$ and initial fluid velocity $u_{0}$
as some measures concentrated on smooth curves and surfaces (manifolds) or
at points.

In view of the inclusions above, our initial data class is larger than that
of Kozono-Miura-Sugiyama \cite{Kozono1}. Moreover, the force field $f$ is
assumed to belong to the Morrey space $\mathcal{M}_{N_{1}}^{N}\left( \mathbb{%
R}^{N}\right) $, where $N_{1}$ can be taken equal to $N\frac{r_{1}}{r}$ with
$r_{1}$ less and close to $r$ and $p/p_{1}=r/r_{1}.$ Thus, in view of (\ref%
{imer-aux-1}), our class of forces $f$ is larger than that of \cite{Kozono1}%
. For $N=3$ and $r\in \lbrack 1,6)$, there exist indexes $p$, $p_{1}$, $q$
and $q_{1}$ satisfying Assumption \ref{Assumption} such that
\begin{equation*}
\dot{B}_{r,1}^{-2+\frac{3}{r}}\hookrightarrow \mathcal{N}_{{q},{q_{1}},{%
\infty }}^{\frac{3}{q}-2},\,\,\,\,\dot{B}_{r,1}^{\frac{3}{r}}\hookrightarrow
\dot{B}_{\infty ,1}^{0}\hookrightarrow L^{\infty },\,\,\,\,\dot{B}_{r,1}^{-1+%
\frac{3}{r}}\hookrightarrow \mathcal{N}_{{p},{p_{1}},{\infty }}^{\frac{3}{p}%
-1},
\end{equation*}%
where $\dot{B}_{p,r}^{s}$ stands for homogeneous Besov spaces ($\dot{B}%
_{p,r}^{s}=\mathcal{N}_{{p},{p},{r}}^{s}$). Then, in the case of (\ref{NSC})
without chemical attractant ($v=0$), our initial data class is larger than
those of \cite{ChoeLkha,ZhaoZhou}. It is worth pointing out that
Besov-Morrey spaces were introduced in \cite{Kozono1994} (see also \cite%
{Mazzucato}) in order to study Navier-Stokes equations.

The mild solutions are obtained by means of a contraction argument in a
time-dependent critical space defined in (\ref{X-space}). Under additional
conditions of homogeneity on the initial data $n_{0}$, $c_{0}$, $v_{0}$, $%
u_{0}$ and the external force $f$, we can ensure that the solution obtained
in Theorem \ref{TeoremaPrincipal} is self-similar when $\gamma =0$. Finally,
we show that solutions are asymptotically stable under small initial
perturbations, as the time goes to infinity. As a byproduct, we obtain a
class of asymptotically self-similar solutions when $\gamma =0.$

This paper is organized as follows. In Section \ref{sec2}, we recall the
definitions of Morrey and Besov-Morrey spaces and present some properties
about these spaces. Our results on well-posedness and asymptotic behavior of
solutions are stated in Section \ref{sec3}. In Section \ref{sec4}, we obtain
the needed linear and nonlinear estimates and prove our results.

\section{Preliminaries}

\label{sec2}

This section is devoted to some preliminaries about Morrey and Besov-Morrey
spaces. For further details about these spaces, see \cite%
{Kato,Kozono1994,Mazzucato}.

\begin{defin}
For $1\leq p_{1}\leq p<\infty $, the Morrey space $\mathcal{M}_{p_{1}}^{p}=%
\mathcal{M}_{p_{1}}^{p}(\mathbb{R}^{N})$ is defined as the set of all
measurable functions $u$ such that
\begin{equation}
\Vert u\Vert _{\mathcal{M}_{p_{1}}^{p}}=\sup_{x_{0}\in \mathbb{R}%
^{N}}\sup_{R>0}R^{\frac{N}{p}-\frac{N}{p_{1}}}\Vert u\Vert
_{L^{p_{1}}(D(x_{0},R))}<\infty ,  \label{morrey}
\end{equation}%
where $D(x_{0},R)$ denotes the closed ball in $\mathbb{R}^{N}$ with center $%
x_{0}$ and radius $R$.
\end{defin}

The space $\mathcal{M}_{p_{1}}^{p}$ endowed with $\Vert \cdot \Vert _{%
\mathcal{M}_{p_{1}}^{p}}$ is a Banach space. In the case $p_{1}=1,$ $%
\mathcal{M}_{1}^{p}$ is a space of signed Radon measures and $\Vert u\Vert
_{L^{1}(D(x_{0},R))}$ is meant as the total variation of the measure $u$ in
the ball $D(x_{0},R)$. For $1<p<\infty $ we have that $\mathcal{M}%
_{p}^{p}=L^{p}$ and $\mathcal{M}_{1}^{1}=\mathcal{M}$ where $\mathcal{M}$
stands for the space of signed Radon measures with finite total variation.
In the case $p=p_{1}=\infty ,$ we consider $\mathcal{M}_{\infty }^{\infty
}=L^{\infty }$.

Next we recall Hölder inequality and heat semigroup estimates in the
framework of Morrey spaces.

\begin{lema}
(Hölder inequality) Let $1\leq p_{1}\leq p\leq \infty $, $1\leq q_{1}\leq
q\leq \infty $ and $1\leq r_{1}\leq r\leq \infty $. If $\frac{1}{r}=\frac{1}{%
p}+\frac{1}{q}$ and $\frac{1}{r_{1}}=\frac{1}{p_{1}}+\frac{1}{q_{1}},$ then
\begin{equation}
\Vert fg\Vert _{\mathcal{M}_{r_{1}}^{r}}\leq \Vert f\Vert _{\mathcal{M}%
_{p_{1}}^{p}}\Vert g\Vert _{\mathcal{M}_{q_{1}}^{q}},  \label{HolderMorrey}
\end{equation}%
for all $f\in \mathcal{M}_{p_{1}}^{p}$ and $g\in \mathcal{M}_{q_{1}}^{q}$.
\end{lema}

Let $\{e^{t\Delta }\}_{t\geq 0}$ denote the heat semigroup. We have the
following estimates in the framework of Morrey spaces.

\begin{lema}
Let $1\leq p_{1}\leq p<\infty $ and $1\leq q_{1}\leq q<\infty $. If $p\geq q$
and $\frac{p}{p_{1}}\geq \frac{q}{q_{1}},$ then there exists a universal
constant $C>0$ such that
\begin{eqnarray}
\Vert e^{t\Delta }f\Vert _{\mathcal{M}_{p_{1}}^{p}} &\leq &C\,t^{-\frac{N}{2}%
(\frac{1}{q}-\frac{1}{p})}\,\Vert f\Vert _{\mathcal{M}_{q_{1}}^{q}}
\label{SemCalorMorrey} \\
\Vert \partial _{x}e^{t\Delta }f\Vert _{\mathcal{M}_{p_{1}}^{p}} &\leq
&C\,t^{-\frac{N}{2}(\frac{1}{q}-\frac{1}{p})-\frac{1}{2}}\,\Vert f\Vert _{%
\mathcal{M}_{q_{1}}^{q}}.  \label{DerivadaSemCalorMorrey}
\end{eqnarray}%
Furthermore, for $1\leq q_{1}\leq q<\infty ,$ it holds that
\begin{eqnarray}
\Vert e^{t\Delta }f\Vert _{L^{\infty }} &\leq &C\,t^{-\frac{N}{2q}}\,\Vert
f\Vert _{\mathcal{M}_{q_{1}}^{q}}  \label{SemCalorMorrey1} \\
\Vert \partial _{x}e^{t\Delta }f\Vert _{L^{\infty }} &\leq &C\,t^{-\frac{N}{%
2q}-\frac{1}{2}}\,\Vert f\Vert _{\mathcal{M}_{q_{1}}^{q}},
\label{DerivadaSemCalorMorrey1}
\end{eqnarray}%
where $C>0$ is a universal constant.
\end{lema}

Let us denote by $\mathcal{S}$ and $\mathcal{S^{\prime }}$ the Schwartz
class and the space of tempered distributions, respectively. For $u\in
\mathcal{S^{\prime }}$, we denote the Fourier transform of $u$ by $\widehat{u%
}$ and its inverse by $u^{\vee }$. Let $\chi (z)$ be a $C^{\infty }$%
-function on $[0,\infty )$ such that $0\leq \chi (z)\leq 1$, $\chi (z)\equiv
1$ for $z\leq 3/2$ and supp $\chi \subset \lbrack 0,5/3)$.Then, for all $%
j\in \mathbb{Z}$, put $\varphi _{j}(\xi )=\chi (2^{-j}|\xi |)-\chi
(2^{1-j}|\xi |)$. It follows that $\varphi _{j}(\xi )\in C_{0}^{\infty }(%
\mathbb{R}^{N})$ and we have the dyadic decomposition
\begin{equation*}
\sum_{j=-\infty }^{\infty }\varphi _{j}(\xi )=1,\,\,\mbox{for all }\xi \neq
0.
\end{equation*}

\begin{defin}
The homogeneous Besov-Morrey space $\mathcal{N}_{{p},{p_{1}},{r}}^{s}=%
\mathcal{N}_{{p},{p_{1}},{r}}^{s}(\mathbb{R}^{N})$ is the set of all $u\in
\mathcal{S^{\prime }}/\mathcal{P}$ such that $\varphi _{j}^{\vee }\ast u\in
\mathcal{M}_{p_{1}}^{p}$ for all $j$, and
\begin{equation}
\Vert u\Vert _{\mathcal{N}_{{p},{p_{1}},{r}}^{s}}=\left\{
\begin{array}{ll}
\left( \displaystyle\sum_{j\in \mathbb{Z}}\left( 2^{sj}\Vert \varphi
_{j}^{\vee }\ast u\Vert _{\mathcal{M}_{p_{1}}^{p}}\right) ^{r}\right) ^{%
\frac{1}{r}}<\infty , & \,\,\,\mbox{for }1\leq p_{1}\leq p\leq \infty
,\,1\leq r<\infty ,\,s\in \mathbb{R}, \\
\displaystyle\sup_{j\in \mathbb{Z}}\left( 2^{sj}\Vert \varphi _{j}^{\vee
}\ast u\Vert _{\mathcal{M}_{p_{1}}^{p}}\right) <\infty , & \,\,\,\mbox{for }%
1\leq p_{1}\leq p\leq \infty ,\,r=\infty ,\,s\in \mathbb{R},%
\end{array}%
\right.  \label{BesovMorrey}
\end{equation}%
where $\mathcal{P}$ denotes the set of polynomials with $N$ variables.
\end{defin}

The space $\mathcal{N}_{{p},{p_{1}},{r}}^{s}$ is a Banach space with the
norm $\Vert \cdot \Vert _{\mathcal{N}_{{p},{p_{1}},{r}}^{s}}$. For all $s\in
\mathbb{R},$ $1\leq p_{2}\leq p_{1}\leq p<\infty $ and $1\leq r\leq \bar{r}%
\leq \infty ,$ we have the continuous inclusions (see \cite{Kozono1994})
\begin{equation}
\mathcal{N}_{{p},{p_{1}},{r}}^{s}\hookrightarrow \mathcal{N}_{{p},{p_{2}},%
\bar{r}}^{s}\,\,\,\,\mbox{and}\,\,\,\,\mathcal{N}_{{p},{p_{1}},{r}%
}^{s}\hookrightarrow \mathcal{N}_{{\frac{p}{\theta }},{\frac{p_{1}}{\theta }}%
,{r}}^{s-\frac{N(1-\theta )}{p}},\,\,\,\mbox{for all}\,\,\,\,\theta \in
(0,1).  \label{imersoes}
\end{equation}%
Furthermore, if $s<0$ we have the following equivalence of norms (see \cite%
{Mazzucato})
\begin{equation}
\Vert u\Vert _{\mathcal{N}_{{p},{p_{1}},{\infty }}^{s}}\cong \sup_{t>0}\,t^{-%
\frac{s}{2}}\Vert e^{t\Delta }u\Vert _{\mathcal{M}_{p_{1}}^{p}}.
\label{EquiNormBesov-morrey}
\end{equation}

\section{Results}

\label{sec3}

In this section we present our global well-posedness and asymptotic behavior
results for the system \eqref{NSC}. Before exposing them, we perform a
scaling analysis for finding the suitable functional setting.

For $\gamma \neq 0$, the system \eqref{NSC} has no scaling relation.
However, we can consider for (\ref{NSC}) the scaling of the case $\gamma =0.$
Assume temporarily that $\gamma =0$ and $f$ is a homogeneous distribution of
degree $-1.$ If $(n,c,v,u,p)$ is a classical solution for (\ref{NSC}), then
so does $(n_{\lambda },c_{\lambda },v_{\lambda },u_{\lambda },p_{\lambda })$
with initial data $(\lambda ^{2}n_{0}(\lambda x),c_{0}(\lambda
x),v_{0}(\lambda x),\lambda \,u_{0}(\lambda x)),$ for each $\lambda >0$,
where $n_{\lambda }(x,t)\coloneqq\lambda ^{2}\,n(\lambda x,\lambda ^{2}t)$, $%
c_{\lambda }(x,t)\coloneqq c(\lambda x,\lambda ^{2}t)$, $v_{\lambda }(x,t)%
\coloneqq v(\lambda x,\lambda ^{2}t)$, $u_{\lambda }(x,t)\coloneqq\lambda
\,u(\lambda x,\lambda ^{2}t)$ and $p_{\lambda }(x,t)\coloneqq\lambda
^{2}\,p(\lambda x,\lambda ^{2}t)$. This leads us to consider the scaling
maps (\ref{scal1}) and (\ref{scal2}).

Solutions invariant by (\ref{scal1}) are named self-similar ones. Then, for $%
(n,c,v,u,p)$ to be a self-similar solution, it is necessary that the initial
data $n_{0},$ $c_{0},$ $v_{0},$ $u_{0}$ and force $f$ are homogeneous
functions of degree $-2,0,0,-1$ and $-1$, respectively. Motivated by the
above scaling analysis, we consider the following critical initial data
class
\begin{eqnarray}
&&n_{0}\in \mathcal{N}_{{q},{q_{1}},{\infty }}^{\frac{N}{q}-2}\left( \mathbb{%
R}^{N}\right) ,\,\,\,\,c_{0}\in L^{\infty }\left( \mathbb{R}^{N}\right)
\text{ with }\nabla c_{0}\in \mathcal{N}_{{r},{r_{1}},{\infty }}^{\frac{N}{r}%
-1}\left( \mathbb{R}^{N}\right) ,  \label{dadosiniciais} \\
&&v_{0}\in \mathcal{S^{\prime }}/\mathcal{P}\text{ with }\nabla v_{0}\in
\mathcal{N}_{{r},{r_{1}},{\infty }}^{\frac{N}{r}-1}\left( \mathbb{R}%
^{N}\right) ,\text{ and }u_{0}\in \mathcal{N}_{{p},{p_{1}},{\infty }}^{\frac{%
N}{p}-1}\left( \mathbb{R}^{N}\right) ,  \notag
\end{eqnarray}%
and force $f\in \mathcal{M}_{N_{1}}^{N}\left( \mathbb{R}^{N}\right) $, where
the exponents $p$, $p_{1}$, $q$, $q_{1}$, $r$, $r_{1}$ and $N_{1}$ are as in
the following assumption.

\begin{Assumption}
\label{Assumption} Assume that $N\geq 2$ and $\gamma \geq 0$. For $N\geq 3$,
suppose that the exponents $p$, $q$ and $r$ satisfy either $(i)$, $(ii)$ or $%
(iii)$ where
\begin{equation*}
\begin{array}{llll}
(i) & \frac{N}{2}<q<N, & N<p<\frac{Nq}{N-q}, & N<r<\frac{Nq}{N-q}; \\
(ii) & q=N, & N<p<\infty , & N<r<\infty ; \\
(iii) & N<q<2N, & N<p<\frac{Nq}{q-N}, & q\leq r<\frac{Nq}{q-N}. \\
&  &  &
\end{array}%
\end{equation*}%
In the case $N=2$ we assume that $p$, $q$ and $r$ satisfy the condition $%
(iii)$ above. Moreover, suppose also that $p_{1}$, $q_{1}$, $r_{1}$ and $%
N_{1}$ satisfy the following conditions
\begin{equation*}
\begin{array}{lllll}
(A) & 1\leq p_{1}\leq p, & \hspace{-2.4cm}1\leq q_{1}\leq q, & 1\leq
r_{1}\leq r, & 1\leq N_{1}\leq N; \\
&  &  &  &  \\
(B) & \displaystyle\frac{1}{p_{1}}+\frac{1}{q_{1}}\leq 1, & \hspace{-2.4cm}%
\displaystyle\frac{1}{r_{1}}+\frac{1}{q_{1}}\leq 1, & \displaystyle\frac{1}{%
p_{1}}+\frac{1}{r_{1}}\leq 1, & \displaystyle\frac{1}{N_{1}}+\frac{1}{q_{1}}%
\leq 1; \\
&  &  &  &  \\
(C) & \displaystyle\frac{p}{p_{1}}\leq \frac{q}{q_{1}}=\frac{r}{r_{1}}; &  &
&  \\
&  &  &  &  \\
(D) & \displaystyle p_{1}\left( \frac{1}{N_{1}}+\frac{1}{q_{1}}\right) \leq
p\left( \frac{1}{N}+\frac{1}{q}\right) . &  &  &  \\
&  &  &  &
\end{array}%
\end{equation*}
\end{Assumption}

\begin{obs}
It is always possible to find indexes $p_{1},q_{1},r_{1}$ and $N_{1}$
sufficiently close to $p,q,r$ and $N$, respectively, satisfying either $(i)$%
, $(ii)$ or $(iii)$, and such that $(A),$ $(B),$ $(C)$ and $(D)$ hold true.
In other words, Assumption \ref{Assumption} is not empty.
\end{obs}

Let $\mathcal{Z}$ be a Banach space continuously included in $\mathcal{S}%
^{\prime }$ and denote by $BC_{w}\left( (0,\infty );\mathcal{Z}\right) $ the
class of bounded functions from $(0,\infty )$ to $\mathcal{Z}$ that are
weakly time continuous in the sense of $\mathcal{S}^{\prime }$. We define
the functional spaces
\begin{eqnarray}
X_{1} &\coloneqq&\left\{ n:t^{-\frac{N}{2q}+1}\,n\in BC_{w}\left( (0,\infty
);\mathcal{M}_{q_{1}}^{q}\right) \right\} ,  \label{X1} \\
X_{2} &\coloneqq&\left\{ c:c\in BC_{w}\left( (0,\infty );L^{\infty }\right)
\,\mbox{ with }t^{-\frac{N}{2r}+\frac{1}{2}}\,\nabla c\in BC_{w}\left(
(0,\infty );\mathcal{M}_{r_{1}}^{r}\right) \right\} ,  \label{X2} \\
X_{3} &\coloneqq&\left\{ v:v(\cdot ,t)\in \mathcal{S^{\prime }}/\mathcal{P}%
\text{ for }t>0\text{ and }t^{-\frac{N}{2r}+\frac{1}{2}}\,\nabla v\in
BC_{w}\left( (0,\infty );\mathcal{M}_{r_{1}}^{r}\right) \right\} ,
\label{X3} \\
X_{4} &\coloneqq&\left\{ u:t^{-\frac{N}{2p}+\frac{1}{2}}\,u\in BC_{w}\left(
(0,\infty );\mathcal{M}_{p_{1}}^{p}\right) \right\} ,  \label{X4}
\end{eqnarray}%
which are Banach spaces endowed with the respective norms
\begin{equation*}
\begin{array}{ll}
& \Vert n\Vert _{X_{1}}\coloneqq\displaystyle\sup_{t>0}\,t^{-\frac{N}{2q}%
+1}\,\Vert n(t)\Vert _{\mathcal{M}_{q_{1}}^{q}},\vspace{0.15cm} \\
& \Vert c\Vert _{X_{2}}\coloneqq\displaystyle\sup_{t>0}\,\Vert c(t)\Vert
_{L^{\infty }}+\displaystyle\sup_{t>0}\,t^{-\frac{N}{2r}+\frac{1}{2}}\,\Vert
\nabla c(t)\Vert _{\mathcal{M}_{r_{1}}^{r}},\vspace{0.15cm} \\
& \Vert v\Vert _{X_{3}}\coloneqq\displaystyle\sup_{t>0}\,t^{-\frac{N}{2r}+%
\frac{1}{2}}\,\Vert \nabla v(t)\Vert _{\mathcal{M}_{r_{1}}^{r}}, \\
& \Vert u\Vert _{X_{4}}\coloneqq\displaystyle\sup_{t>0}\,t^{-\frac{N}{2p}+%
\frac{1}{2}}\,\Vert u(t)\Vert _{\mathcal{M}_{p_{1}}^{p}}\vspace{0.15cm}.%
\end{array}%
\end{equation*}

Next, let us introduce the spaces $\mathcal{X}$ and $\mathcal{I}$ as
\begin{equation}
\mathcal{X}\coloneqq\left\{ (n,c,v,u):n\in X_{1},\,c\in X_{2},\,v\in
X_{3},\,u\in X_{4}\right\}  \label{X-space}
\end{equation}%
with the norm
\begin{equation*}
\left\Vert (n,c,v,u)\right\Vert _{\mathcal{X}}\coloneqq\Vert n\Vert
_{X_{1}}+\,\Vert c\Vert _{X_{2}}+\,\Vert v\Vert _{X_{3}}+\,\Vert u\Vert
_{X_{4}},
\end{equation*}%
and
\begin{equation*}
\mathcal{I}\coloneqq\left\{ (n_{0},c_{0},v_{0},u_{0}):n_{0},\,c_{0},v_{0}%
\text{ and }u_{0}\text{ are as in (\ref{dadosiniciais})}\right\}
\end{equation*}%
with the norm
\begin{equation*}
\left\Vert (n_{0},c_{0},v_{0},u_{0})\right\Vert _{\mathcal{I}}\coloneqq\Vert
n_{0}\Vert _{\mathcal{N}_{{q},{q_{1}},{\infty }}^{\frac{N}{q}-2}}+\,\Vert
c_{0}\Vert _{L^{\infty }}+\,\Vert \nabla c_{0}\Vert _{\mathcal{N}_{{r},{r_{1}%
},{\infty }}^{\frac{N}{r}-1}}+\,\Vert \nabla v_{0}\Vert _{\mathcal{N}_{{r},{%
r_{1}},{\infty }}^{\frac{N}{r}-1}}+\,\Vert u_{0}\Vert _{\mathcal{N}_{{p},{%
p_{1}},{\infty }}^{\frac{N}{p}-1}}.
\end{equation*}%
Note that $\mathcal{X}$ and $\mathcal{I}$ are Banach spaces equipped with
the norms $\Vert \cdot \Vert _{\mathcal{X}}$ and $\Vert \cdot \Vert _{%
\mathcal{I}}$, respectively.

Let $\mathbb{P}=I+\mathcal{R}\otimes \mathcal{R}$ stand for the
Leray-Helmholtz projection onto the spaces of solenoidal vector fields,
where $I$ is the identity and $\mathcal{R}=(\mathcal{R}_{1},\mathcal{R}%
_{1},\ldots ,\mathcal{R}_{N})$ is a vector of operators whose components are
the Riesz transforms $\mathcal{R}_{j}$. Applying $\mathbb{P}$ on the fourth
equation in \eqref{NSC} and using Duhamel's principle, system \eqref{NSC}
can be formally converted to the following integral formulation
\begin{equation}
\left\{
\begin{array}{lll}
n(t) & = & e^{t\Delta }n_{0}-\displaystyle\int_{0}^{t}e^{(t-\tau )\Delta
}(u\cdot \nabla n)(\tau )\,d\tau -\displaystyle\int_{0}^{t}\nabla \cdot
e^{(t-\tau )\Delta }(n\nabla c+n\nabla v)(\tau )\,d\tau , \\
c(t) & = & e^{t\Delta }c_{0}-\displaystyle\int_{0}^{t}e^{(t-\tau )\Delta
}(u\cdot \nabla c+nc)(\tau )\,d\tau , \\
v(t) & = & e^{-\gamma t}e^{t\Delta }v_{0}-\displaystyle\int_{0}^{t}e^{-%
\gamma (t-\tau )}e^{(t-\tau )\Delta }(u\cdot \nabla v-n)(\tau )\,d\tau , \\
u(t) & = & e^{t\Delta }u_{0}-\displaystyle\int_{0}^{t}e^{(t-\tau )\Delta }%
\mathbb{P}(u\cdot \nabla u)\,d\tau -\displaystyle\int_{0}^{t}e^{(t-\tau
)\Delta }\mathbb{P}(nf)(\tau )\,d\tau .%
\end{array}%
\right.  \label{EI}
\end{equation}

A 4-tuple $(n,c,v,u)$ satisfying \eqref{EI} is called a \textit{mild solution%
} of \eqref{NSC}. In what follows, we state our main result.

\begin{teo}
\label{TeoremaPrincipal} Let $N\geq 2$, and let the exponents $p$, $p_{1}$, $%
q$, $q_{1}$, $r$, $r_{1}$ and $N_{1}$ be as in Assumption \ref{Assumption}.
Suppose that the initial data $(n_{0},c_{0},v_{0},u_{0})\in \mathcal{I}$ and
the external force $f\in \mathcal{M}_{N_{1}}^{N}(\mathbb{R}^{N})$. There
exist positive constants $\varepsilon ,$ $\delta $ ($\delta =C\varepsilon $)
and $K_{1}$ such that the system \eqref{EI} has a unique global mild
solution $(n,c,v,u)\in \mathcal{X}$ satisfying $\Vert (n,c,v,u)\Vert _{%
\mathcal{X}}\leq 2K_{1}\varepsilon $ provided that $\Vert
(n_{0},c_{0},v_{0},u_{0})\Vert _{\mathcal{I}}\leq \delta $. Moreover, the
data-solution map is locally Lipschitz continuous.
\end{teo}

\begin{obs}
Let the constants $C_{i}$'s$,$ $i=1,...,7,$ be as in Lemma \ref%
{LemaEstimativaBilinear} and $\alpha ,\beta $ as in Lemma \ref%
{LemaEstimativaLinear}. The constant $\varepsilon $ in Theorem \ref%
{TeoremaPrincipal} can be chosen so that $0<\varepsilon <\frac{1}{4K_{1}K_{2}%
},$ where $K_{1}$ and $K_{2}$ depend on $C_{i},\alpha ,\beta $ (see %
\eqref{K1eK2}). We also point out that the mild solution $%
(n,c,v,u)\rightharpoonup (n_{0},c_{0},v_{0},u_{0})$ in the sense of
distributions, as $t\rightarrow 0^{+}$.
\end{obs}

Since the space $\mathcal{X}$ is critical with respect to the scaling of the
case $\gamma =0,$ we can obtain self-similar solutions by assuming the right
homogeneity on the data and force.

\begin{cor}
(Self-similar solution) \label{self-similar1}Let $N\geq 3$ and $\gamma =0$.
Assume that $(n_{0},c_{0},v_{0},u_{0})$ and $f$ are as in Theorem \ref%
{TeoremaPrincipal}. Suppose that $n_{0},c_{0},v_{0},u_{0}$ and $f$ are
homogeneous functions with degree $-2,0,0,-1$ and $-1$, respectively. Then,
the solution $(n,c,v,u)$ obtained through Theorem \ref{TeoremaPrincipal} is
self-similar, that is, for every $\lambda >0$ we have that
\begin{equation*}
n(x,t)=\lambda ^{2}n(\lambda x,\lambda ^{2}t),\,\,c(x,t)=c(\lambda x,\lambda
^{2}t),\,\,v(x,t)=v(\lambda x,\lambda ^{2}t)\,\text{ and }\,u(x,t)=\,\lambda
u(\lambda x,\lambda ^{2}t).
\end{equation*}
\end{cor}

Now we present an asymptotic stability result for solutions of system $%
\eqref{NSC}$.

\begin{teo}
\label{casymptotic behavior} Under the hypotheses of Theorem \ref%
{TeoremaPrincipal}. Assume that $(n,c,v,u)$ and $(\tilde{n},\tilde{c},\tilde{%
v},\tilde{u})$ are two solutions given by Theorem \ref{TeoremaPrincipal}
corresponding to the initial data $(n_{0},c_{0},v_{0},u_{0})$ and $(\tilde{n}%
_{0},\tilde{c}_{0},\tilde{v}_{0},\tilde{u}_{0})$, respectively. We have that%
\begin{eqnarray}
\lim_{t\rightarrow \infty }t^{-\frac{N}{2q}+1}\Vert n(\cdot ,t)-\tilde{n}%
(\cdot ,t)\Vert _{\mathcal{M}_{q_{1}}^{q}} &=&\lim_{t\rightarrow \infty
}\Vert c(\cdot ,t)-\tilde{c}(\cdot ,t)\Vert _{L^{\infty
}}=\lim_{t\rightarrow \infty }t^{-\frac{N}{2r}+\frac{1}{2}}\Vert \nabla
(c(\cdot ,t)-\tilde{c}(\cdot ,t))\Vert _{\mathcal{M}_{r_{1}}^{r}}=  \notag \\
\lim_{t\rightarrow \infty }t^{-\frac{N}{2r}+\frac{1}{2}}\Vert \nabla
(v(\cdot ,t)-\tilde{v}(\cdot ,t))\Vert _{\mathcal{M}_{r_{1}}^{r}}
&=&\lim_{t\rightarrow \infty }t^{-\frac{N}{2p}+\frac{1}{2}}\Vert u(\cdot ,t)-%
\tilde{u}(\cdot ,t)\Vert _{\mathcal{M}_{p_{1}}^{p}}=0.  \label{volta}
\end{eqnarray}
if only if
\begin{eqnarray}
&&\lim_{t\rightarrow \infty }\Big(t^{-\frac{N}{2q}+1}\Vert e^{t\Delta
}(n_{0}-\tilde{n}_{0})\Vert _{\mathcal{M}_{q_{1}}^{q}}+\Vert e^{t\Delta
}(c_{0}-\tilde{c}_{0})\Vert _{L^{\infty }}+t^{-\frac{N}{2r}+\frac{1}{2}%
}\Vert \nabla e^{t\Delta }(c_{0}-\tilde{c}_{0})\Vert _{\mathcal{M}%
_{r_{1}}^{r}}+  \notag \\
&&t^{-\frac{N}{2r}+\frac{1}{2}}\Vert \nabla e^{-\gamma t}e^{t\Delta }(v_{0}-%
\tilde{v}_{0})\Vert _{\mathcal{M}_{r_{1}}^{r}}+t^{-\frac{N}{2p}+\frac{1}{2}%
}\Vert e^{t\Delta }(u_{0}-\tilde{u}_{0})\Vert _{\mathcal{M}_{p_{1}}^{p}}\Big)%
=0,  \label{ida}
\end{eqnarray}
\end{teo}

\begin{obs}
(Asymptotically self-similar solutions) In the case $\gamma =0,$ Theorem \ref%
{casymptotic behavior} together with Corollary \ref{self-similar1} provide a
class of solutions asymptotically self-similar at infinity. Indeed, taking
the initial data $(\tilde{n}_{0},\tilde{c}_{0},\tilde{v}_{0},\tilde{u}%
_{0})=(n_{0},c_{0},v_{0},u_{0})+(\varphi _{1},\varphi _{2},\varphi
_{3},\varphi _{4})$ with $\varphi _{i}\in C_{0}^{\infty }$ and $%
n_{0},c_{0},v_{0},u_{0}$ and $f$ as in Corollary \ref{self-similar1}, we
have that the corresponding solution $(\tilde{n},\tilde{c},\tilde{v},\tilde{u%
})$ is attracted to the self-similar solution $(n,c,v,u)$ in the sense of %
\eqref{volta}.
\end{obs}

\section{Proofs}

\label{sec4}

In this section we present the proofs of the results stated in Section \ref%
{sec3}. First we prove an abstract fixed point lemma which will be useful
for our ends.

\begin{lema}
\label{teoDoPontoFixo} For $1\leq i\leq 4,$ let $X_{i}$ be a Banach space
with the norm $\Vert \cdot \Vert _{X_{i}}$. Consider the Banach space $%
\mathcal{X}=X_{1}\times X_{2}\times X_{3}\times X_{4}$ endowed with the norm
\begin{equation*}
\left\Vert x\right\Vert _{\mathcal{X}}=\left\Vert x_{1}\right\Vert
_{X_{1}}+\left\Vert x_{2}\right\Vert _{X_{2}}+\left\Vert x_{3}\right\Vert
_{X_{3}}+\left\Vert x_{4}\right\Vert _{X_{4}},
\end{equation*}%
where $x=(x_{1},x_{2},x_{3},x_{4})\in \mathcal{X}$. For $1\leq i,j,k\leq 4,$
assume that $B_{i,j}^{k}:X_{i}\times X_{j}\rightarrow X_{k}$ is a continuous
bilinear map, that is, there is a constant $C_{i,j}^{k}>0$ such that
\begin{equation}
\left\Vert B_{i,j}^{k}(x_{i},x_{j})\right\Vert _{X_{k}}\leq
C_{i,j}^{k}\,\left\Vert x_{i}\right\Vert _{X_{i}}\,\left\Vert
x_{j}\right\Vert _{X_{j}},\,\,\,\,\,\,\,\,\,\,\mbox{for
all }(x_{i},x_{j})\in X_{i}\times X_{j}.  \label{bili-abstract-1}
\end{equation}%
Assume also that $L_{3}:X_{1}\rightarrow X_{3}$ and $L_{4}:X_{1}\rightarrow
X_{4}$ are continuous linear maps such that $\left\Vert L_{3}\right\Vert
_{X_{1}\rightarrow X_{3}}\leq \alpha $ and $\left\Vert L_{4}\right\Vert
_{X_{1}\rightarrow X_{4}}\leq \beta $. Set the constants
\begin{equation*}
K_{1}\coloneqq1+\alpha +\beta \,\,\,\,\,\,\mbox{and}\,\,\,\,\,\,K_{2}%
\coloneqq(\alpha +\beta
)\sum_{i,j=1}^{4}C_{i,j}^{1}+\sum_{k,i,j=1}^{4}C_{i,j}^{k},
\end{equation*}%
and let $0<\varepsilon <\displaystyle\frac{1}{4K_{1}K_{2}}$ and $\mathcal{B}%
_{\varepsilon }=\left\{ x\in \mathcal{X}:\Vert x\Vert _{\mathcal{X}}\leq
2\,K_{1}\,\varepsilon \right\} $. If $\Vert y\Vert _{\mathcal{X}}\leq
\varepsilon $ then there exists a unique solution $x\in \mathcal{B}%
_{\varepsilon }$ for the equation $x=y+B(x)$, where $%
y=(y_{1},y_{2},y_{3},y_{4}),$ $B(x)=\left(
B_{1}(x),B_{2}(x),B_{3}(x),B_{4}(x)\right) $ and
\begin{eqnarray*}
B_{1}(x) &=&\sum_{i,j=1}^{4}B_{i,j}^{1}(x_{i},x_{j}), \\
B_{2}(x) &=&\sum_{i,j=1}^{4}B_{i,j}^{2}(x_{i},x_{j}), \\
B_{3}(x) &=&\sum_{i,j=1}^{4}B_{i,j}^{3}(x_{i},x_{j})+\left( L_{3}\circ
(y_{1}+B_{1})\right) (x), \\
B_{4}(x) &=&\sum_{i,j=1}^{4}B_{i,j}^{4}(x_{i},x_{j})+\left( L_{4}\circ
(y_{1}+B_{1})\right) (x).
\end{eqnarray*}
\end{lema}

\textbf{Proof.} For all $x\in \mathcal{X}$, it follows from (\ref%
{bili-abstract-1}) that
\begin{eqnarray}
\left\Vert B_{1}(x)\right\Vert _{X_{1}} &\leq &\sum_{i,j=1}^{4}\left\Vert
B_{i,j}^{1}(x_{i},x_{j})\right\Vert _{X_{1}}  \notag \\
&\leq &\sum_{i,j=1}^{4}C_{i,j}^{1}\,\Vert x_{i}\Vert _{X_{i}}\,\Vert
x_{j}\Vert _{X_{j}}  \notag \\
&\leq &\left( \sum_{i,j=1}^{4}C_{i,j}^{1}\right) \,\Vert x\Vert _{\mathcal{X}%
}^{2}.  \label{estB1}
\end{eqnarray}%
Analogously, we have
\begin{equation*}
\left\Vert B_{2}(x)\right\Vert _{X_{2}}\leq \left(
\sum_{i,j=1}^{4}C_{i,j}^{2}\right) \,\Vert x\Vert _{\mathcal{X}}^{2}.
\end{equation*}%
Next, using (\ref{bili-abstract-1}) and (\ref{estB1}), we estimate $B_{3}$
as follows:
\begin{eqnarray}
\left\Vert B_{3}(x)\right\Vert _{X_{3}} &\leq &\sum_{i,j=1}^{4}\left\Vert
B_{i,j}^{3}(x_{i},x_{j})\right\Vert _{X_{3}}+\left\Vert \left( L_{3}\circ
(y_{1}+B_{1})\right) (x)\right\Vert _{X_{3}}  \notag \\
&\leq &\sum_{i,j=1}^{4}C_{i,j}^{3}\,\Vert x_{i}\Vert _{X_{i}}\,\Vert
x_{j}\Vert _{X_{j}}+\alpha \,\left\Vert (y_{1}+B_{1})(x)\right\Vert _{X_{1}}
\notag \\
&\leq &\left( \sum_{i,j=1}^{4}C_{i,j}^{3}\right) \,\Vert x\Vert _{\mathcal{X}%
}^{2}+\alpha \,\left( \Vert y\Vert _{\mathcal{X}}+\left(
\sum_{i,j=1}^{4}C_{i,j}^{1}\right) \,\Vert x\Vert _{\mathcal{X}}^{2}\right)
\notag \\
&=&\left( \sum_{i,j=1}^{4}C_{i,j}^{3}+\alpha
\,\sum_{i,j=1}^{4}C_{i,j}^{1}\right) \,\Vert x\Vert _{\mathcal{X}%
}^{2}+\alpha \,\Vert y\Vert _{\mathcal{X}}.  \label{estB3}
\end{eqnarray}%
Similarly, it follows that
\begin{equation}
\left\Vert B_{4}(x)\right\Vert _{X_{4}}\leq \left(
\sum_{i,j=1}^{4}C_{i,j}^{4}+\beta \,\sum_{i,j=1}^{4}C_{i,j}^{1}\right)
\,\Vert x\Vert _{\mathcal{X}}^{2}+\beta \,\Vert y\Vert _{\mathcal{X}}.
\label{estB4}
\end{equation}

Now, consider the mapping $\mathcal{F}:\mathcal{X}\rightarrow \mathcal{X}$
given by $F(x)=y+B(x)$. For $x\in \mathcal{B}_{\varepsilon }$, from (\ref%
{estB1})-(\ref{estB4}) we obtain that
\begin{eqnarray*}
\left\Vert \mathcal{F}(x)\right\Vert _{\mathcal{X}} &\leq &\Vert y\Vert _{%
\mathcal{X}}+\sum_{k=1}^{4}\left\Vert B_{k}(x)\right\Vert _{X_{k}} \\
&\leq &(1+\alpha +\beta )\Vert y\Vert _{\mathcal{X}}+\left( (\alpha +\beta
)\sum_{i,j=1}^{4}C_{i,j}^{1}+\sum_{k,i,j=1}^{4}C_{i,j}^{k}\right) \,\Vert
x\Vert _{\mathcal{X}}^{2} \\
&\leq &K_{1}\,\varepsilon +K_{2}\,4\,K_{1}^{2}\,\varepsilon
^{2}=(1+4\,K_{1}\,K_{2}\,\varepsilon )\,K_{1}\,\varepsilon \leq
2\,K_{1}\,\varepsilon ,
\end{eqnarray*}%
and then $\mathcal{F}\left( \mathcal{B}_{\varepsilon }\right) \subset
\mathcal{B}_{\varepsilon }$. Next, we take $x,z\in \mathcal{B}_{\varepsilon
} $ and estimate
\begin{eqnarray*}
\left\Vert \mathcal{F}(x)-\mathcal{F}(z)\right\Vert _{\mathcal{X}}
&=&\left\Vert B(x)-B(z)\right\Vert _{\mathcal{X}} \\
&=&\sum_{k=1}^{4}\left\Vert B_{k}(x)-B_{k}(z)\right\Vert _{X_{k}} \\
&\leq &\sum_{k=1}^{4}\sum_{i,j=1}^{4}\left\Vert
B_{i,j}^{k}(x_{i}-z_{i},x_{j})+B_{i,j}^{k}(z_{i},x_{j}-z_{j})\right\Vert
_{X_{k}}+\sum_{l=3}^{4}\left\Vert \left( L_{l}\circ
(B_{1}(x)-B_{1}(z))\right) \right\Vert _{X_{l}} \\
&\leq &\sum_{k,i,j=1}^{4}C_{i,j}^{k}\,\left( \Vert x_{i}-z_{i}\Vert
_{X_{i}}\,\Vert x_{j}\Vert _{X_{j}}+\Vert z_{i}\Vert _{X_{i}}\,\Vert
x_{j}-z_{j}\Vert _{X_{j}}\right) +(\alpha +\beta )\left\Vert
B_{1}(x)-B_{1}(z)\right\Vert _{X_{1}} \\
&\leq &\sum_{k,i,j=1}^{4}C_{i,j}^{k}\,\Vert x-z\Vert _{\mathcal{X}}\left(
\Vert x\Vert _{\mathcal{X}}+\Vert z\Vert _{\mathcal{X}}\right) \\
&&+\,\,(\alpha +\beta )\sum_{i,j=1}^{4}C_{i,j}^{1}\,\left( \Vert
x_{i}-z_{i}\Vert _{X_{i}}\,\Vert x_{j}\Vert _{X_{j}}+\Vert z_{i}\Vert
_{X_{i}}\,\Vert x_{j}-z_{j}\Vert _{X_{j}}\right) \\
&\leq &\left( (\alpha +\beta
)\sum_{i,j=1}^{4}C_{i,j}^{1}+\sum_{k,i,j=1}^{4}C_{i,j}^{k}\right) \,\Vert
x-z\Vert _{\mathcal{X}}\left( \Vert x\Vert _{\mathcal{X}}+\Vert z\Vert _{%
\mathcal{X}}\right) \\
&\leq &K_{2}\,4\,K_{1}\,\varepsilon \,\Vert x-z\Vert _{\mathcal{X}}.
\end{eqnarray*}%
Since $4K_{1}K_{2}\varepsilon <1$, $\mathcal{F}$ is a contraction in $%
\mathcal{B}_{\varepsilon }$, and the Banach fixed point theorem concludes
the proof. \begin{flushright}$\blacksquare$\end{flushright}

\begin{obs}
{\label{obs}} Due to the fixed point argument in the proof of Lemma \ref%
{teoDoPontoFixo}, we have that the solution $x$ depends continuously on the
data $y$. More precisely, the data-solution map is Lipschitz continuous from
$\{y\in \mathcal{X};\left\Vert y\right\Vert \leq \varepsilon \}$ to $%
\mathcal{B}_{\varepsilon }.$ In addition, the solution obtained through
Lemma \ref{teoDoPontoFixo} is the limit in $\mathcal{X}$ of the sequence of
iterates $x^{(1)}=y$ and $x^{(m+1)}=\mathcal{F}(x^{(m)})$, $m\geq 1$. This
fact will be useful in the proof of Corollary \ref{self-similar1}.
\end{obs}

Now, for each initial data tuple $(n_{0},c_{0},v_{0},u_{0})$ and force $f,$
we consider $\mathcal{F}(n,c,v,u)=(\mathcal{N},\mathcal{C},\mathcal{V},%
\mathcal{U})$, where
\begin{equation}
\left\{
\begin{array}{lll}
\mathcal{N}(t) & = & e^{t\Delta }n_{0}-\displaystyle\int_{0}^{t}e^{(t-\tau
)\Delta }(u\cdot \nabla n)(\tau )\,d\tau -\displaystyle\int_{0}^{t}\nabla
\cdot e^{(t-\tau )\Delta }(n\nabla c)(\tau )\,d\tau -\displaystyle%
\int_{0}^{t}\nabla \cdot e^{(t-\tau )\Delta }(n\nabla v)(\tau )\,d\tau , \\
&  &  \\
& =: & e^{t\Delta }n_{0}+B_{4,1}^{1}(t)+B_{1,2}^{1}(t)+B_{1,3}^{1}(t), \\
&  &  \\
\mathcal{C}(t) & = & e^{t\Delta }c_{0}-\displaystyle\int_{0}^{t}e^{(t-\tau
)\Delta }(u\cdot \nabla c)(\tau )\,d\tau -\displaystyle\int_{0}^{t}e^{(t-%
\tau )\Delta }(nc)(\tau )\,d\tau , \\
&  &  \\
& =: & e^{t\Delta }c_{0}+B_{4,2}^{2}(t)+B_{1,2}^{2}(t), \\
&  &  \\
\mathcal{V}(t) & = & e^{-\gamma t}e^{t\Delta }v_{0}-\displaystyle%
\int_{0}^{t}e^{-\gamma (t-\tau )}e^{(t-\tau )\Delta }(u\cdot \nabla v)(\tau
)\,d\tau +\displaystyle\int_{0}^{t}e^{-\gamma (t-\tau )}e^{(t-\tau )\Delta
}n(\tau )\,d\tau , \\
&  &  \\
& =: & e^{-\gamma t}e^{t\Delta }v_{0}+B_{4,3}^{3}(t)+L_{3}(t), \\
&  &  \\
\mathcal{U}(t) & = & e^{t\Delta }u_{0}-\displaystyle\int_{0}^{t}e^{(t-\tau
)\Delta }\mathbb{P}(u\cdot \nabla u)\,d\tau -\displaystyle%
\int_{0}^{t}e^{(t-\tau )\Delta }\mathbb{P}(nf)(\tau )\,d\tau , \\
&  &  \\
& =: & e^{t\Delta
}u_{0}+B_{4,4}^{4}(t)+L_{4}(t),\,\,\,\,\,\,\,\,\,\,0<t<\infty . \\
&  &
\end{array}%
\right.  \label{NCVU}
\end{equation}

\subsection{Estimates for the bilinear terms in \eqref{NCVU}}

\begin{lema}
\label{LemaEstimativaBilinear} Under the hypotheses of Theorem \ref%
{TeoremaPrincipal}. There exist positive constants $%
C_{1},C_{2},C_{3},C_{4},C_{5},C_{6},C_{7}$ such that
\begin{eqnarray}
\left\Vert B_{4,1}^{1}(u,n)\right\Vert _{X_{1}} &\leq &C_{1}\,\Vert u\Vert
_{X_{4}}\,\Vert n\Vert _{X_{1}},  \label{B1_4,1} \\
\left\Vert B_{1,2}^{1}(n,c)\right\Vert _{X_{1}} &\leq &C_{2}\,\Vert n\Vert
_{X_{1}}\,\Vert c\Vert _{X_{2}},  \label{B1_1,2} \\
\left\Vert B_{1,3}^{1}(n,v)\right\Vert _{X_{1}} &\leq &C_{3}\,\Vert n\Vert
_{X_{1}}\,\Vert v\Vert _{X_{3}},  \label{B1_1,3} \\
\left\Vert B_{4,2}^{2}(u,c)\right\Vert _{X_{2}} &\leq &C_{4}\,\Vert u\Vert
_{X_{4}}\,\Vert c\Vert _{X_{2}},  \label{B2_4,2} \\
\left\Vert B_{1,2}^{2}(n,c)\right\Vert _{X_{2}} &\leq &C_{5}\,\Vert n\Vert
_{X_{1}}\,\Vert c\Vert _{X_{2}},  \label{B2_1,2} \\
\left\Vert B_{4,3}^{3}(u,v)\right\Vert _{X_{3}} &\leq &C_{6}\,\Vert u\Vert
_{X_{4}}\,\Vert v\Vert _{X_{3}},  \label{B3_4,3} \\
\left\Vert B_{4,4}^{4}(u,\tilde{u})\right\Vert _{X_{4}} &\leq &C_{7}\,\Vert
u\Vert _{X_{4}}\Vert \tilde{u}\Vert _{X_{4}},  \label{B4_4,4}
\end{eqnarray}%
for all $n\in X_{1},c\in X_{2},v\in X_{3}$ and $u,\tilde{u}\in X_{4}.$
\end{lema}

\textbf{Proof.} From the conditions $(i)$, $(ii)$ and $(iii)$ in Assumption $%
\ref{Assumption}$, we have that
\begin{equation*}
\frac{1}{2}-\frac{N}{2p}>0,-\frac{1}{2}+\frac{N}{2p}+\frac{N}{2q}>0.
\end{equation*}%
Taking $s_{1}=\frac{p_{1}q_{1}}{p_{1}+q_{1}},$ from $(A)$, $(B)$ and $(C)$
in Assumption \ref{Assumption}, it follows that
\begin{equation*}
1\leq s_{1}\leq \frac{pq}{p+q}\leq q\text{ and }\frac{q}{q_{1}}\geq \frac{pq%
}{p+q}\frac{1}{s_{1}},
\end{equation*}%
and hence we can estimate
\begin{eqnarray}
\left\Vert B_{4,1}^{1}(u,n)(t)\right\Vert _{\mathcal{M}_{q_{1}}^{q}}
&=&\left\Vert \int_{0}^{t}e^{(t-\tau )\Delta }(u\cdot \nabla n)(\tau
)\,d\tau \right\Vert _{\mathcal{M}_{q_{1}}^{q}}  \notag \\
&\leq &\int_{0}^{t}\left\Vert \nabla \cdot e^{(t-\tau )\Delta }(un)(\tau
)\right\Vert _{\mathcal{M}_{q_{1}}^{q}}\,d\tau  \notag \\
&\leq &C\,\int_{0}^{t}(t-\tau )^{-\frac{N}{2}(\frac{1}{q}+\frac{1}{p}-\frac{1%
}{q})-\frac{1}{2}}\,\left\Vert (un)(\tau )\right\Vert _{\mathcal{M}_{s_{1}}^{%
\frac{pq}{p+q}}}\,d\tau \,\,\,(\mbox{by }\eqref{DerivadaSemCalorMorrey})
\notag \\
&\leq &C\,\int_{0}^{t}(t-\tau )^{-\frac{N}{2p}-\frac{1}{2}}\,\left\Vert
u(\tau )\right\Vert _{\mathcal{M}_{p_{1}}^{p}}\left\Vert n(\tau )\right\Vert
_{\mathcal{M}_{q_{1}}^{q}}\,d\tau \,\,\,(\mbox{by }\eqref{HolderMorrey})
\notag \\
&\leq &C\,\int_{0}^{t}(t-\tau )^{-\frac{N}{2p}-\frac{1}{2}}\,\tau ^{\frac{N}{%
2p}-\frac{1}{2}}\,\tau ^{\frac{N}{2q}-1}\,d\tau \,\Vert u\Vert
_{X_{4}}\,\Vert n\Vert _{X_{1}}  \notag \\
&=&C\,t^{\frac{N}{2q}-1}b\left( \frac{1}{2}-\frac{N}{2p},-\frac{1}{2}+\frac{N%
}{2p}+\frac{N}{2q}\right) \,\Vert u\Vert _{X_{4}}\,\Vert n\Vert _{X_{1}}
\notag \\
&=&C_{1}\,t^{\frac{N}{2q}-1}\,\Vert u\Vert _{X_{4}}\,\Vert n\Vert _{X_{1}},
\label{estB^1_4,1}
\end{eqnarray}%
for all $t>0,$ where $C_{1}=C(N,p,p_{1},q,q_{1})$ and $b(\cdot ,\cdot )$
denotes the beta function.

Taking $s_{2}=\frac{r_{1}q_{1}}{r_{1}+q_{1}}$, we have that
\begin{equation*}
\frac{1}{2}-\frac{N}{2r}>0,-\frac{1}{2}+\frac{N}{2q}+\frac{N}{2r}>0,1\leq
s_{2}\leq \frac{rq}{r+q}\leq q\text{ \ and }\frac{q}{q_{1}}\geq \frac{rq}{r+q%
}\frac{1}{s_{2}},
\end{equation*}%
and then
\begin{eqnarray}
\left\Vert B_{1,2}^{1}(n,c)(t)\right\Vert _{\mathcal{M}_{q_{1}}^{q}}
&=&\left\Vert \int_{0}^{t}\nabla \cdot e^{(t-\tau )\Delta }(n\nabla c)(\tau
)\,d\tau \right\Vert _{\mathcal{M}_{q_{1}}^{q}}  \notag \\
&\leq &\int_{0}^{t}\left\Vert \nabla \cdot e^{(t-\tau )\Delta }(n\nabla
c)(\tau )\right\Vert _{\mathcal{M}_{q_{1}}^{q}}\,d\tau  \notag \\
&\leq &C\,\int_{0}^{t}(t-\tau )^{-\frac{N}{2}(\frac{1}{q}+\frac{1}{r}-\frac{1%
}{q})-\frac{1}{2}}\,\left\Vert (n\nabla c)(\tau )\right\Vert _{\mathcal{M}%
_{s_{2}}^{\frac{rq}{r+q}}}\,d\tau \,\,\,(\mbox{by }%
\eqref{DerivadaSemCalorMorrey})  \notag \\
&\leq &C\,\int_{0}^{t}(t-\tau )^{-\frac{N}{2r}-\frac{1}{2}}\,\left\Vert
n(\tau )\right\Vert _{\mathcal{M}_{q_{1}}^{q}}\,\left\Vert \nabla c(\tau
)\right\Vert _{\mathcal{M}_{r_{1}}^{r}}\,d\tau \,\,\,(\mbox{by }%
\eqref{HolderMorrey})  \notag \\
&\leq &C\,\int_{0}^{t}(t-\tau )^{-\frac{N}{2r}-\frac{1}{2}}\,\tau ^{\frac{N}{%
2q}-1}\,\tau ^{\frac{N}{2r}-\frac{1}{2}}\,d\tau \,\Vert n\Vert
_{X_{1}}\,\Vert c\Vert _{X_{2}}  \notag \\
&=&C\,t^{\frac{N}{2q}-1}\,\Vert n\Vert _{X_{1}}\,\Vert c\Vert
_{X_{2}}b\left( \frac{1}{2}-\frac{N}{2r},-\frac{1}{2}+\frac{N}{2q}+\frac{N}{%
2r}\right)  \notag \\
&=&C_{2}\,t^{\frac{N}{2q}-1}\,\Vert n\Vert _{X_{1}}\,\Vert c\Vert _{X_{2}},
\label{estB^1_1,2}
\end{eqnarray}%
for all $t>0,$ where $C_{2}=C(N,q,q_{1},r,r_{1})$. Similarly,
\begin{eqnarray}
\left\Vert B_{1,3}^{1}(n,v)(t)\right\Vert _{\mathcal{M}_{q_{1}}^{q}} &\leq
&C\,t^{\frac{N}{2q}-1}\,\Vert n\Vert _{X_{1}}\,\Vert v\Vert
_{X_{3}}\,b\left( \frac{1}{2}-\frac{N}{2r},-\frac{1}{2}+\frac{N}{2q}+\frac{N%
}{2r}\right)  \notag \\
&=&C_{3}\,t^{\frac{N}{2q}-1}\,\Vert n\Vert _{X_{1}}\,\Vert v\Vert _{X_{3}},
\label{estB^1_1,3}
\end{eqnarray}%
for all $t>0,$ where $C_{3}=C_{3}(N,q,q_{1},r,r_{1})$. Thus, from (\ref%
{estB^1_4,1})-(\ref{estB^1_1,3}), we obtain the inequalities \eqref{B1_4,1}, %
\eqref{B1_1,2} and \eqref{B1_1,3}.

Now, from $(i)$, $(ii)$ and $(iii)$ in Assumption \ref{Assumption}, we have
that
\begin{equation*}
\frac{1}{2}-\frac{N}{2p}>0,1-\frac{N}{2q}>0\text{ and }\frac{1}{2}-\frac{N}{%
2q}+\frac{N}{2r}>0.
\end{equation*}%
Taking $s_{3}=\frac{p_{1}r_{1}}{p_{1}+r_{1}},$ the conditions $(A)$, $(B)$
and $(C)$ in Assumption \ref{Assumption} gives that
\begin{equation*}
1\leq s_{3}\leq \frac{pr}{p+r}\leq r\text{ and }\frac{r}{r_{1}}\geq \frac{pr%
}{p+r}\frac{1}{s_{3}}.
\end{equation*}%
Hence, we have the following estimates for $B_{4,2}^{2}$ and $\nabla
B_{4,2}^{2}$:
\begin{eqnarray}
\left\Vert B_{4,2}^{2}(u,c)(t)\right\Vert _{L^{\infty }} &=&\left\Vert
\int_{0}^{t}e^{(t-\tau )\Delta }(u\cdot \nabla c)(\tau )\,d\tau \right\Vert
_{L^{\infty }}  \notag \\
&\leq &\int_{0}^{t}\left\Vert \nabla \cdot e^{(t-\tau )\Delta }(uc)(\tau
)\right\Vert _{L^{\infty }}\,d\tau  \notag \\
&\leq &C\,\int_{0}^{t}(t-\tau )^{-\frac{N}{2p}-\frac{1}{2}}\,\left\Vert
(uc)(\tau )\right\Vert _{\mathcal{M}_{p_{1}}^{p}}\,d\tau \,\,\,(\mbox{by }%
\eqref{DerivadaSemCalorMorrey1})  \notag \\
&\leq &C\,\int_{0}^{t}(t-\tau )^{-\frac{N}{2p}-\frac{1}{2}}\,\left\Vert
u(\tau )\right\Vert _{\mathcal{M}_{p_{1}}^{p}}\,\left\Vert c(\tau
)\right\Vert _{L^{\infty }}\,d\tau  \notag \\
&\leq &C\,\int_{0}^{t}(t-\tau )^{-\frac{N}{2p}-\frac{1}{2}}\,\tau ^{\frac{N}{%
2p}-\frac{1}{2}}\,d\tau \,\Vert u\Vert _{X_{4}}\,\Vert c\Vert _{X_{2}}
\notag \\
&=&C\,\Vert u\Vert _{X_{4}}\,\Vert c\Vert _{X_{2}}b\left( \frac{1}{2}-\frac{N%
}{2p},\frac{1}{2}+\frac{N}{2p}\right)  \notag \\
&=&C_{4,1}\,\Vert u\Vert _{X_{4}}\,\Vert c\Vert _{X_{2}},  \label{estB^2_4,2}
\end{eqnarray}%
and
\begin{eqnarray}
\left\Vert \nabla B_{4,2}^{2}(u,c)(t)\right\Vert _{\mathcal{M}_{r_{1}}^{r}}
&=&\left\Vert \nabla \int_{0}^{t}e^{(t-\tau )\Delta }(u\cdot \nabla c)(\tau
)\,d\tau \right\Vert _{\mathcal{M}_{r_{1}}^{r}}  \notag \\
&\leq &\int_{0}^{t}\left\Vert \nabla e^{(t-\tau )\Delta }(u\cdot \nabla
c)(\tau )\right\Vert _{\mathcal{M}_{r_{1}}^{r}}\,d\tau  \notag \\
&\leq &C\,\int_{0}^{t}(t-\tau )^{-\frac{N}{2}(\frac{1}{p}+\frac{1}{r}-\frac{1%
}{r})-\frac{1}{2}}\,\left\Vert (u\cdot \nabla c)(\tau )\right\Vert _{%
\mathcal{M}_{s_{3}}^{\frac{pr}{p+r}}}\,d\tau \,\,\,(\mbox{by }%
\eqref{DerivadaSemCalorMorrey})  \notag \\
&\leq &C\,\int_{0}^{t}(t-\tau )^{-\frac{N}{2p}-\frac{1}{2}}\,\left\Vert
u(\tau )\right\Vert _{\mathcal{M}_{p_{1}}^{p}}\,\left\Vert \nabla c(\tau
)\right\Vert _{\mathcal{M}_{r_{1}}^{r}}\,d\tau \,\,\,(\mbox{by }%
\eqref{HolderMorrey})  \notag \\
&\leq &C\,\int_{0}^{t}(t-\tau )^{-\frac{N}{2p}-\frac{1}{2}}\,\tau ^{\frac{N}{%
2p}-\frac{1}{2}}\,\tau ^{\frac{N}{2r}-\frac{1}{2}}\,d\tau \,\Vert u\Vert
_{X_{4}}\,\Vert c\Vert _{X_{2}}  \notag \\
&\leq &C\,t^{\frac{N}{2r}-\frac{1}{2}}b\left( \frac{1}{2}-\frac{N}{2p},\frac{%
N}{2p}+\frac{N}{2r}\right) \,\Vert u\Vert _{X_{4}}\,\Vert c\Vert _{X_{2}}
\notag \\
&=&C_{4,2}\,t^{\frac{N}{2r}-\frac{1}{2}}\,\Vert u\Vert _{X_{4}}\,\Vert
c\Vert _{X_{2}}.  \label{estDeltaB^2_4,2}
\end{eqnarray}%
In turn, we can estimate $B_{1,2}^{2}$ and $\nabla B_{1,2}^{2}$ as follows:%
\begin{eqnarray}
\left\Vert B_{1,2}^{2}(n,c)(t)\right\Vert _{L^{\infty }} &=&\left\Vert
\int_{0}^{t}e^{(t-\tau )\Delta }(nc)(\tau )\,d\tau \right\Vert _{L^{\infty }}
\notag \\
&\leq &\int_{0}^{t}\left\Vert e^{(t-\tau )\Delta }(nc)(\tau )\right\Vert
_{L^{\infty }}\,d\tau  \notag \\
&\leq &C\,\int_{0}^{t}(t-\tau )^{-\frac{N}{2q}}\,\left\Vert (nc)(\tau
)\right\Vert _{\mathcal{M}_{q_{1}}^{q}}\,d\tau \,\,\,(\mbox{by }%
\eqref{SemCalorMorrey1})  \notag \\
&\leq &C\,\int_{0}^{t}(t-\tau )^{-\frac{N}{2q}}\,\left\Vert n(\tau
)\right\Vert _{\mathcal{M}_{q_{1}}^{q}}\,\left\Vert c(\tau )\right\Vert
_{L^{\infty }}\,d\tau  \notag \\
&\leq &C\,\int_{0}^{t}(t-\tau )^{-\frac{N}{2q}}\,\tau ^{\frac{N}{2q}%
-1}\,d\tau \,\Vert n\Vert _{X_{1}}\,\Vert c\Vert _{X_{2}}  \notag \\
&=&C\,\Vert n\Vert _{X_{1}}\,\Vert c\Vert _{X_{2}}b\left( 1-\frac{N}{2q},%
\frac{N}{2q}\right)  \notag \\
&=&C_{5,1}\,\Vert n\Vert _{X_{1}}\,\Vert c\Vert _{X_{2}},  \label{estB^2_1,2}
\end{eqnarray}%
and
\begin{eqnarray}
\left\Vert \nabla B_{1,2}^{2}(n,c)(t)\right\Vert _{\mathcal{M}_{r_{1}}^{r}}
&=&\left\Vert \nabla \int_{0}^{t}e^{(t-\tau )\Delta }(nc)(\tau )\,d\tau
\right\Vert _{\mathcal{M}_{r_{1}}^{r}}  \notag \\
&\leq &\int_{0}^{t}\left\Vert \nabla e^{(t-\tau )\Delta }(nc)(\tau
)\right\Vert _{\mathcal{M}_{r_{1}}^{r}}\,d\tau  \notag \\
&\leq &C\,\int_{0}^{t}(t-\tau )^{-\frac{N}{2}(\frac{1}{q}-\frac{1}{r})-\frac{%
1}{2}}\,\left\Vert (nc)(\tau )\right\Vert _{\mathcal{M}_{q_{1}}^{q}}\,d\tau
\,\,\,(\mbox{by }\eqref{DerivadaSemCalorMorrey})  \notag \\
&\leq &C\,\int_{0}^{t}(t-\tau )^{-\frac{N}{2q}+\frac{N}{2r}-\frac{1}{2}%
}\,\left\Vert n(\tau )\right\Vert _{\mathcal{M}_{q_{1}}^{q}}\,\left\Vert
c(\tau )\right\Vert _{L^{\infty }}\,d\tau  \notag \\
&\leq &C\,\int_{0}^{t}(t-\tau )^{-\frac{N}{2q}+\frac{N}{2r}-\frac{1}{2}%
}\,\tau ^{\frac{N}{2q}-1}\,d\tau \,\Vert n\Vert _{X_{1}}\,\Vert c\Vert
_{X_{2}}  \notag \\
&\leq &C\,t^{\frac{N}{2r}-\frac{1}{2}}b\left( \frac{1}{2}-\frac{N}{2q}+\frac{%
N}{2r},\frac{N}{2q}\right) \,\Vert n\Vert _{X_{1}}\,\Vert c\Vert _{X_{2}}
\notag \\
&=&C_{5,2}\,t^{\frac{N}{2r}-\frac{1}{2}}\,\Vert n\Vert _{X_{1}}\,\Vert
c\Vert _{X_{2}},  \label{estDeltaB^2_1,2}
\end{eqnarray}%
for all $t>0,$ where $C_{4,1}=C_{4,1}(N,p,p_{1})$, $%
C_{4,2}=C_{4,2}(N,p,p_{1},r,r_{1})$, $C_{5,1}=C_{5,1}(N,q,q_{1})$ and $%
C_{5,2}=C_{5,2}(N,q,q_{1},r,r_{1})$. Taking $C_{4}=C_{4,1}+C_{4,2}$ and $%
C_{5}=C_{5,1}+C_{5,2},$ estimates \eqref{B2_4,2} and \eqref{B2_1,2} follow
from \eqref{estB^2_4,2}-\eqref{estDeltaB^2_1,2}.

Proceeding similarly to (\ref{estB^2_4,2}), we can estimate $\nabla
B_{4,3}^{3}$ in $\mathcal{M}_{r_{1}}^{r}$ as
\begin{eqnarray}
\left\Vert \nabla B_{4,3}^{3}(u,v)(t)\right\Vert _{\mathcal{M}_{r_{1}}^{r}}
&=&\left\Vert \nabla \int_{0}^{t}e^{-\gamma (t-\tau )}\,e^{(t-\tau )\Delta
}\,(u\cdot \nabla v)(\tau )\,d\tau \right\Vert _{\mathcal{M}_{r_{1}}^{r}}
\notag \\
&\leq &C\,t^{\frac{N}{2r}-\frac{1}{2}}b\left( \frac{1}{2}-\frac{N}{2p},\frac{%
N}{2p}+\frac{N}{2r}\right) \,\Vert u\Vert _{X_{4}}\Vert v\Vert _{X_{3}}
\notag \\
&=&C_{6}\,t^{\frac{N}{2r}-\frac{1}{2}}\,\Vert u\Vert _{X_{4}}\Vert v\Vert
_{X_{3}},  \label{estB^3_4,3}
\end{eqnarray}%
for all $t>0,$ where $C_{6}=C_{6}(N,p,p_{1},r,r_{1})$, which gives %
\eqref{B3_4,3}.

Finally, since the projection operator $\mathbb{P}$ is bounded in $\mathcal{M%
}_{p_{1}}^{p}$, we have that
\begin{eqnarray}
\left\Vert B_{4,4}^{4}(u,\tilde{u})(t)\right\Vert _{\mathcal{M}_{p_{1}}^{p}}
&=&\left\Vert \int_{0}^{t}e^{(t-\tau )\Delta }\,\mathbb{P}(u\cdot \nabla
\tilde{u})\,d\tau \right\Vert _{\mathcal{M}_{p_{1}}^{p}}  \notag \\
&\leq &\int_{0}^{t}\left\Vert \mathbb{P}\,\nabla \cdot e^{(t-\tau )\Delta
}(u\otimes \tilde{u})(\tau )\right\Vert _{\mathcal{M}_{p_{1}}^{p}}\,d\tau
\notag \\
&\leq &C\,\int_{0}^{t}\left\Vert \nabla \cdot e^{(t-\tau )\Delta
}\,(u\otimes \tilde{u})(\tau )\right\Vert _{\mathcal{M}_{p_{1}}^{p}}\,d\tau
\notag \\
&\leq &C\,\int_{0}^{t}(t-\tau )^{-\frac{N}{2}(\frac{2}{p}-\frac{1}{p})-\frac{%
1}{2}}\,\left\Vert (u\otimes \tilde{u})(\tau )\right\Vert _{\mathcal{M}_{%
\frac{p_{1}}{2}}^{\frac{p}{2}}}\,d\tau \,\,\,(\mbox{by }%
\eqref{DerivadaSemCalorMorrey})  \notag \\
&\leq &C\,\int_{0}^{t}(t-\tau )^{-\frac{N}{2p}-\frac{1}{2}}\,\left\Vert
u(\tau )\right\Vert _{\mathcal{M}_{p_{1}}^{p}}\left\Vert \tilde{u}(\tau
)\right\Vert _{\mathcal{M}_{p_{1}}^{p}}\,d\tau \,\,\,(\mbox{by }%
\eqref{HolderMorrey})  \notag \\
&\leq &C\,\int_{0}^{t}(t-\tau )^{-\frac{N}{2p}-\frac{1}{2}}\,\tau ^{\frac{N}{%
p}-1}\,d\tau \,\Vert u\Vert _{X^{4}}\Vert \tilde{u}\Vert _{X^{4}}  \notag \\
&\leq &C\,t^{\frac{N}{2p}-\frac{1}{2}}b\left( \frac{1}{2}-\frac{N}{2p},\frac{%
N}{p}\right) \,\Vert u\Vert _{X^{4}}\Vert \tilde{u}\Vert _{X^{4}}  \notag \\
&=&C_{7}\,t^{\frac{N}{2p}-\frac{1}{2}}\,\Vert u\Vert _{X^{4}}\Vert \tilde{u}%
\Vert _{X^{4}},  \label{estB^4_4,4}
\end{eqnarray}%
for all $t>0,$ where $C_{7}=C_{7}(N,p,p_{1})$, and then we obtain %
\eqref{B4_4,4}. \begin{flushright}$\blacksquare$\end{flushright}

\subsection{Estimates for the linear terms in (\protect\ref{NCVU})}

\begin{lema}
\label{LemaEstimativaLinear} Under the hypotheses of Theorem \ref%
{TeoremaPrincipal}. There exist constants $\alpha ,\beta >0$ such that
\begin{eqnarray}
\left\Vert L_{3}(n)\right\Vert _{X_{3}} &\leq &\alpha \,\Vert n\Vert
_{X_{1}},  \label{L3} \\
\left\Vert L_{4}(n)\right\Vert _{X_{4}} &\leq &\beta \,\Vert n\Vert _{X_{1}},
\label{L4}
\end{eqnarray}%
for all $n\in X_{1}$.
\end{lema}

\textbf{Proof.} Using \eqref{DerivadaSemCalorMorrey}, we can estimate
\begin{eqnarray}
\left\Vert \nabla L_{3}(n)(t)\right\Vert _{\mathcal{M}_{r_{1}}^{r}}
&=&\left\Vert \nabla \int_{0}^{t}e^{-\gamma (t-\tau )}\,e^{(t-\tau )\Delta
}\,n(\tau )\,d\tau \right\Vert _{\mathcal{M}_{r_{1}}^{r}}  \notag \\
&\leq &\int_{0}^{t}\left\Vert \nabla e^{(t-\tau )\Delta }n(\tau )\right\Vert
_{\mathcal{M}_{r_{1}}^{r}}\,d\tau  \notag \\
&\leq &C\,\int_{0}^{t}(t-\tau )^{-\frac{N}{2}(\frac{1}{q}-\frac{1}{r})-\frac{%
1}{2}}\,\left\Vert n(\tau )\right\Vert _{\mathcal{M}_{q_{1}}^{q}}\,d\tau
\notag \\
&\leq &C\,t^{\frac{N}{2r}-\frac{1}{2}}\,b\left( \frac{1}{2}-\frac{N}{2q}+%
\frac{N}{2r},\frac{N}{2q}\right) \,\Vert n\Vert _{X_{1}}  \notag \\
&=&\alpha \,t^{\frac{N}{2r}-\frac{1}{2}}\,\Vert n\Vert _{X_{1}},
\label{estL_3}
\end{eqnarray}%
for all $t>0,$ where $\alpha =\alpha (N,q,q_{1},r,r_{1})$, which gives (\ref%
{L3}).

Now, considering $s_{4}=\frac{N_{1}q_{1}}{N_{1}+q_{1}},$ from $(A)$, $(B)$
and $(D)$ in Assumption \ref{Assumption}, we have that
\begin{equation*}
1\leq s_{4}\leq \frac{Nq}{N+q}\leq p\text{ and }\frac{p}{p_{1}}\geq \frac{Nq%
}{N+q}\frac{1}{s_{4}}.
\end{equation*}%
Thus, $L_{4}(n)$ can be estimated as follows:
\begin{eqnarray}
\left\Vert L_{4}(n)(t)\right\Vert _{\mathcal{M}_{p_{1}}^{p}} &=&\left\Vert
\int_{0}^{t}e^{(t-\tau )\Delta }\mathbb{P}(nf)(\tau )\,d\tau \right\Vert _{%
\mathcal{M}_{p_{1}}^{p}}  \notag \\
&\leq &\int_{0}^{t}\left\Vert \mathbb{P}\,e^{(t-\tau )\Delta }(nf)(\tau
)\right\Vert _{\mathcal{M}_{p_{1}}^{p}}\,d\tau  \notag \\
&\leq &C\,\int_{0}^{t}(t-\tau )^{-\frac{N}{2}(\frac{1}{N}+\frac{1}{q}-\frac{1%
}{p})}\,\left\Vert nf(\tau )\right\Vert _{\mathcal{M}_{s_{4}}^{\frac{Nq}{N+q}%
}}\,d\tau \,\,\,(\mbox{by }\eqref{SemCalorMorrey})  \notag \\
&\leq &C\,\int_{0}^{t}(t-\tau )^{-\frac{N}{2q}+\frac{N}{2p}-\frac{1}{2}%
}\,\left\Vert f\right\Vert _{\mathcal{M}_{N_{1}}^{N}}\,\left\Vert n(\tau
)\right\Vert _{\mathcal{M}_{q_{1}}^{q}}\,d\tau \,\,\,(\mbox{by }%
\eqref{HolderMorrey})  \notag \\
&\leq &C\,t^{\frac{N}{2p}-\frac{1}{2}}\,\Vert f\Vert _{\mathcal{M}%
_{N_{1}}^{N}}\,\Vert n\Vert _{X_{1}}\,b\left( \frac{1}{2}+\frac{N}{2p}-\frac{%
N}{2q},\frac{N}{2q}\right)  \notag \\
&=&\beta \,t^{\frac{N}{2p}-\frac{1}{2}}\,\Vert n\Vert _{X_{1}},
\label{estL_4}
\end{eqnarray}%
for all $t>0,$ where $\beta =\beta (N,N_{1},p,p_{1},q,q_{1},f)$, as
requested. \begin{flushright}$\blacksquare$\end{flushright}

\subsection{Proof of Theorem\protect\ref{TeoremaPrincipal}}

Consider $X_{1},X_{2},X_{3}$ and $X_{4}$ as in \eqref{X1}-\eqref{X4} and let
$y=\left( e^{t\Delta }n_{0},e^{t\Delta }c_{0},e^{-\gamma t}e^{t\Delta
}v_{0},e^{t\Delta }u_{0}\right) $. For $\mathcal{X}=X_{1}\times X_{2}\times
X_{3}\times X_{4}$ and $x=(n,c,v,u)\in \mathcal{X},$ we denote
\begin{eqnarray}
B_{1}(x) &\coloneqq&B_{4,1}^{1}(u,n)+B_{1,2}^{1}(n,c)+B_{1,3}^{1}(n,v),
\label{B1} \\
B_{2}(x) &\coloneqq&B_{4,2}^{2}(u,c)+B_{1,2}^{2}(n,c),  \label{B2} \\
B_{3}(x) &\coloneqq&B_{4,3}^{3}(u,v)+L_{3}\circ \left( e^{t\Delta
}n_{0}+B_{1}\right) (x),  \label{B3} \\
B_{4}(x) &\coloneqq&B_{4,4}^{4}(u,u)+L_{4}\circ \left( e^{t\Delta
}n_{0}+B_{1}\right) (x).  \label{B4}
\end{eqnarray}%
From Lemma \ref{LemaEstimativaBilinear}, the operators $B_{i,j}^{k}$ in %
\eqref{B1}-\eqref{B4} are continuous bilinear maps. Also, from Lemma \ref%
{LemaEstimativaLinear}, $L_{3}$ and $L_{4}$ are continuous linear maps.
Moreover, it is not difficult to see that all of them are time-weakly
continuous at $t>0.$

Next, we set
\begin{equation}
K_{1}=1+\alpha +\beta \text{ and }K_{2}=(\alpha +\beta
)(C_{1}+C_{2}+C_{3})+\sum_{i=1}^{7}C_{i}.  \label{K1eK2}
\end{equation}%
In view of equivalence \eqref{EquiNormBesov-morrey}, we have that
\begin{eqnarray}
\Vert y\Vert _{\mathcal{X}} &=&\Vert e^{t\Delta }n_{0}\Vert _{X_{1}}+\Vert
e^{t\Delta }c_{0}\Vert _{X_{2}}+\Vert e^{-\gamma t}e^{t\Delta }v_{0}\Vert
_{X_{3}}+\Vert e^{t\Delta }u_{0}\Vert _{X_{4}}  \notag \\
&=&\sup_{t>0}\,t^{-\frac{N}{2q}+1}\,\Vert e^{t\Delta }n_{0}\Vert _{\mathcal{M%
}_{q_{1}}^{q}}+\sup_{t>0}\,\Vert e^{t\Delta }c_{0}\Vert _{L^{\infty
}}+\sup_{t>0}\,t^{-\frac{N}{2r}+\frac{1}{2}}\,\Vert \nabla e^{t\Delta
}c_{0}\Vert _{\mathcal{M}_{r_{1}}^{r}}  \notag \\
&&+\sup_{t>0}\,t^{-\frac{N}{2r}+\frac{1}{2}}\,\Vert \nabla e^{-\gamma
t}e^{t\Delta }v_{0}\Vert _{\mathcal{M}_{r_{1}}^{r}}+\sup_{t>0}\,t^{-\frac{N}{%
2p}+\frac{1}{2}}\,\Vert e^{t\Delta }u_{0}\Vert _{\mathcal{M}_{p_{1}}^{p}}
\notag \\
&\leq &C_{0}\left( \Vert n_{0}\Vert _{\mathcal{N}_{{q},{q_{1}},{\infty }}^{%
\frac{N}{q}-2}}+\,\Vert c_{0}\Vert _{L^{\infty }}+\,\Vert \nabla c_{0}\Vert
_{\mathcal{N}_{{r},{r_{1}},{\infty }}^{\frac{N}{r}-1}}+\,\Vert \nabla
v_{0}\Vert _{\mathcal{N}_{{r},{r_{1}},{\infty }}^{\frac{N}{r}-1}}+\,\Vert
u_{0}\Vert _{\mathcal{N}_{{p},{p_{1}},{\infty }}^{\frac{N}{p}-1}}\right)
\notag \\
&=&C_{0}\Vert (n_{0},c_{0},v_{0},u_{0})\Vert _{\mathcal{I}}\leq \varepsilon
\label{aux-estimate-initial-1}
\end{eqnarray}%
provided that $\Vert (n_{0},c_{0},v_{0},u_{0})\Vert _{\mathcal{I}}\leq
\delta =\frac{\varepsilon }{C_{0}}$. If $0<\varepsilon <\frac{1}{4K_{1}K_{2}}
$, then Lemma \ref{teoDoPontoFixo} implies that there exists a unique
solution $(n,c,v,u)\in \mathcal{X}$ of (\ref{EI}) such that $\Vert
(n,c,v,u)\Vert _{\mathcal{X}}\leq 2K_{1}\varepsilon $. The continuity of the
data-solution map follows from Remark \ref{obs} and estimate (\ref%
{aux-estimate-initial-1}).\begin{flushright}$\blacksquare$\end{flushright}

\subsection{Proof of Corollary \protect\ref{self-similar1}}

Since we use a fixed point argument to prove Theorem \ref{TeoremaPrincipal},
the solution $(n,c,v,u)$ is the limit in the space $\mathcal{X}$ of the
following Picard sequence (see Remark \ref{obs}):
\begin{equation*}
(n^{(1)},c^{(1)},v^{(1)},u^{(1)})=(e^{t\Delta }n_{0},e^{t\Delta
}c_{0},e^{-\gamma t}e^{t\Delta }v_{0},e^{t\Delta }u_{0})
\end{equation*}%
and
\begin{equation*}
(n^{(m+1)},c^{(m+1)},v^{(m+1)},u^{(m+1)})=(n^{(1)},c^{(1)},v^{(1)},u^{(1)})+%
\mathcal{F}(n^{(m)},c^{(m)},v^{(m)},u^{(m)}),\text{ for }m\in \mathbb{N}.
\end{equation*}%
In other words,
\begin{equation*}
\left\{
\begin{array}{lll}
n^{(m+1)} & = & e^{t\Delta }n_{0}-\displaystyle\int_{0}^{t}e^{(t-\tau
)\Delta }(u^{(m)}\cdot \nabla n^{(m)})(\tau )\,d\tau -\displaystyle%
\int_{0}^{t}\nabla \cdot e^{(t-\tau )\Delta }(n^{(m)}\nabla
c^{(m)}+n^{(m)}\nabla v^{(m)})(\tau )\,d\tau , \\
&  &  \\
c^{(m+1)} & = & e^{t\Delta }c_{0}-\displaystyle\int_{0}^{t}e^{(t-\tau
)\Delta }(u^{(m)}\cdot \nabla c^{(m)}+n^{(m)}c^{(m)})(\tau )\,d\tau , \\
&  &  \\
v^{(m+1)} & = & e^{-\gamma t}e^{t\Delta }v_{0}-\displaystyle%
\int_{0}^{t}e^{-\gamma (t-\tau )}e^{(t-\tau )\Delta }(u^{(m)}\cdot \nabla
v^{(m)}-n^{(m)})(\tau )\,d\tau , \\
&  &  \\
u^{(m+1)} & = & e^{t\Delta }u_{0}-\displaystyle\int_{0}^{t}e^{(t-\tau
)\Delta }\mathbb{P}(u^{(m)}\cdot \nabla u^{(m)})\,d\tau -\displaystyle%
\int_{0}^{t}e^{(t-\tau )\Delta }\mathbb{P}(n^{(m)}f)(\tau )\,d\tau . \\
&  &
\end{array}%
\right.
\end{equation*}

By hypotheses we have that $n_{0}$, $c_{0}$, $v_{0}$, $u_{0}$ and $f$ are
homogeneous functions of degree $-2$, $0$, $0$, $-1$ and $-1$, respectively.
Then, through a simple computation we can verify that $%
(n^{(1)},c^{(1)},v^{(1)},u^{(1)})$ is invariant by (\ref{scal1}), that is,
\begin{eqnarray}
n^{(1)}(x,t) &=&\lambda ^{2}\,n^{(1)}(\lambda x,\lambda
^{2}t),\,\,c^{(1)}(x,t)=c^{(1)}(\lambda x,\lambda ^{2}t),\,\,
\label{invariancia} \\
v^{(1)}(x,t) &=&v^{(1)}(\lambda x,\lambda ^{2}t)\text{ and }%
u^{(1)}(x,t)=\lambda \,u^{(1)}(\lambda x,\lambda ^{2}t).  \notag
\end{eqnarray}%
By means of an induction argument, we can check that $%
(n^{(m)},c^{(m)},v^{(m)},u^{(m)})$ also satisfies the scaling property (\ref%
{invariancia}), for all $m$. Since $(n,c,v,u)$ is the limit in $\mathcal{X}$
of the sequence $\left\{ (n^{(m)},c^{(m)},v^{(m)},u^{(m)})\right\} _{m\in
\mathbb{N}}$ and the norm $\Vert \cdot \Vert _{\mathcal{X}}$ is scaling
invariant, we obtain that the solution $(n,c,v,u)$ is self-similar.%
\begin{flushright}$\blacksquare$\end{flushright}

\subsection{Proof of Theorem \protect\ref{casymptotic behavior}}

We first show that (\ref{ida}) implies (\ref{volta}). Let $(n,c,v,u)$ and $(%
\tilde{n},\tilde{c},\tilde{v},\tilde{u})$ be two mild solutions given by
Theorem \ref{TeoremaPrincipal} and set
\begin{equation*}
l_{q}=-\frac{N}{2q}+1,\,\,\,\mu _{r}=-\frac{N}{2r}+\frac{1}{2}\text{ and }%
\mu _{p}=-\frac{N}{2p}+\frac{1}{2}.
\end{equation*}%
Estimating the difference $n-\tilde{n}$ in the norm $t^{l_{q}}\,\Vert \cdot
\Vert _{\mathcal{M}_{q_{1}}^{q}}$, we obtain
\begin{eqnarray}
t^{-\frac{N}{2q}+1}\Vert n(t)-\tilde{n}(t)\Vert _{\mathcal{M}_{q_{1}}^{q}}
&\leq &t^{l_{q}}\Vert e^{t\Delta }(n_{0}-\tilde{n}_{0})\Vert _{\mathcal{M}%
_{q_{1}}^{q}}  \notag \\
&&+t^{l_{q}}\displaystyle\int_{0}^{t}\Vert e^{(t-\tau )\Delta }(u\cdot
\nabla n-\tilde{u}\cdot \nabla \tilde{n})(\tau )\Vert _{\mathcal{M}%
_{q_{1}}^{q}}\,d\tau  \notag \\
&&+\,t^{l_{q}}\displaystyle\int_{0}^{t}\Vert \nabla \cdot e^{(t-\tau )\Delta
}(n\nabla c+n\nabla v-\tilde{n}\nabla \tilde{c}-\tilde{n}\nabla \tilde{v}%
)(\tau )\Vert _{\mathcal{M}_{q_{1}}^{q}}\,d\tau  \notag \\
&\coloneqq&t^{l_{q}}\,\Vert e^{t\Delta }(n_{0}-\tilde{n}_{0})\Vert _{%
\mathcal{M}_{q_{1}}^{q}}+J_{1}(t)+J_{2}(t).  \label{n-n}
\end{eqnarray}%
The integral $J_{1}$ is estimated as follows:
\begin{eqnarray}
J_{1}(t) &\leq &\tilde{C_{1}}\,t^{l_{q}}\int_{0}^{t}(t-\tau )^{\mu _{p}-1}%
\big(\Vert (u-\tilde{u})(\tau )\Vert _{\mathcal{M}_{p_{1}}^{p}}\Vert n(\tau
)\Vert _{\mathcal{M}_{q_{1}}^{q}}  \notag \\
\text{ \ \ \ } &&\text{ \ \ \ \ \ \ \ \ }+\Vert \tilde{u}(\tau )\Vert _{%
\mathcal{M}_{p_{1}}^{p}}\Vert (n-\tilde{n})(\tau )\Vert _{\mathcal{M}%
_{q_{1}}^{q}}\big)\,d\tau  \notag \\
&\leq &\tilde{C_{1}}\,t^{l_{q}}\int_{0}^{t}(t-\tau )^{\mu _{p}-1}\,\tau
^{-\mu _{p}-l_{q}}\,\tau ^{\mu _{p}}\Vert (u-\tilde{u})(\tau )\Vert _{%
\mathcal{M}_{p_{1}}^{p}}\Vert n\Vert _{X_{1}}\,d\tau  \notag \\
&&+\,\tilde{C_{1}}\,t^{l_{q}}\int_{0}^{t}(t-\tau )^{\mu _{p}-1}\,\tau ^{-\mu
_{p}-l_{q}}\,\tau ^{l_{q}}\Vert \tilde{u}\Vert _{X_{4}}\Vert (n-\tilde{n}%
)(\tau )\Vert _{\mathcal{M}_{q_{1}}^{q}}\,d\tau ,\mbox{ taking }\tau =tz
\notag \\
&=&\tilde{C_{1}}\,\int_{0}^{1}(1-z)^{\mu _{p}-1}\,z^{-\mu
_{p}-l_{q}}\,(tz)^{\mu _{p}}\Vert (u-\tilde{u})(tz)\Vert _{\mathcal{M}%
_{p_{1}}^{p}}\Vert n\Vert _{X_{1}}\,dz  \notag \\
&&+\,\tilde{C_{1}}\,\int_{0}^{1}(1-z)^{\mu _{p}-1}\,z^{-\mu
_{p}-l_{q}}\,(tz)^{l_{q}}\Vert \tilde{u}\Vert _{X_{4}}\Vert (n-\tilde{n}%
)(tz)\Vert _{\mathcal{M}_{q_{1}}^{q}}\,dz.  \label{J1}
\end{eqnarray}

Similarly, from \eqref{estB^1_1,2} and \eqref{estB^1_1,3} we arrive at
\begin{eqnarray}
J_{2}(t) &\leq &\tilde{C_{2}}\,\int_{0}^{1}(1-z)^{\mu _{r}-1}\,z^{-l_{q}-\mu
_{r}}\,\big((tz)^{l_{q}}\,\Vert (n-\tilde{n})(tz)\Vert _{\mathcal{M}%
_{q_{1}}^{q}}\,\Vert c\Vert _{X_{2}}  \notag \\
&&\text{ \ \ \ \ }+(tz)^{\mu _{r}}\,\Vert \nabla (c-\tilde{c})(tz)\Vert _{%
\mathcal{M}_{r_{1}}^{r}}\,\Vert \tilde{n}\Vert _{X_{1}}\,\big)dz  \notag \\
&&+\,\tilde{C_{3}}\int_{0}^{1}(1-z)^{\mu _{r}-1}\,z^{-l_{q}-\mu _{r}}\,\big(%
(tz)^{l_{q}}\,\Vert (n-\tilde{n})(tz)\,\Vert _{\mathcal{M}_{q_{1}}^{q}}\Vert
v\Vert _{X_{3}}  \notag \\
&&\text{ \ \ \ \ \ \ \ \ }+(tz)^{\mu _{r}}\Vert \nabla (v-\tilde{v}%
)(tz)\Vert _{\mathcal{M}_{r_{1}}^{r}}\Vert \tilde{n}\Vert _{X_{1}}\,\big)%
dz.\,\,\,\,\,\,\,\,\,\,\,\,\,\,\,  \label{J2}
\end{eqnarray}

In the following we estimate the differences $c-\tilde{c}$, $v-\tilde{v}$
and $u-\tilde{u}$ in the norms $\Vert \cdot \Vert _{L^{\infty }}+t^{\mu
_{r}}\,\Vert \nabla \cdot \Vert _{\mathcal{M}_{r_{1}}^{r}}$, $t^{\mu
_{r}}\,\Vert \nabla \cdot \Vert _{\mathcal{M}_{r_{1}}^{r}}$ and $t^{\mu
_{p}}\,\Vert \cdot \Vert _{\mathcal{M}_{p_{1}}^{p}}$, respectively. In this
direction, we obtain
\begin{eqnarray}
\Vert (c-\tilde{c})(t)\Vert _{L^{\infty }} &\leq &\Vert e^{t\Delta }(c_{0}-%
\tilde{c}_{0})\Vert _{L^{\infty }}+\int_{0}^{t}\Vert e^{(t-\tau )\Delta
}(u\cdot \nabla c+nc-\tilde{u}\cdot \nabla \tilde{c}-\tilde{n}\tilde{c}%
)(\tau )\Vert _{L^{\infty }}\,d\tau  \notag \\
&\coloneqq&\Vert e^{t\Delta }(c_{0}-\tilde{c}_{0})\Vert _{L^{\infty
}}+J_{3}(t),  \label{c-c(1)}
\end{eqnarray}%
\vspace{-0.6cm}
\begin{eqnarray}
t^{\mu _{r}}\,\Vert \nabla (c-\tilde{c})(t)\Vert _{\mathcal{M}_{r_{1}}^{r}}
&\leq &t^{\mu _{r}}\,\Vert \nabla e^{t\Delta }(c_{0}-\tilde{c}_{0})\Vert _{%
\mathcal{M}_{r_{1}}^{r}}  \notag \\
&&+t^{\mu _{r}}\int_{0}^{t}\Vert \nabla e^{(t-\tau )\Delta }(u\cdot \nabla
c+nc-\tilde{u}\cdot \nabla \tilde{c}-\tilde{n}\tilde{c})(\tau )\Vert _{%
\mathcal{M}_{r_{1}}^{r}}\,d\tau  \notag \\
&\coloneqq&t^{\mu _{r}}\,\Vert \nabla e^{t\Delta }(c_{0}-\tilde{c}_{0})\Vert
_{\mathcal{M}_{r_{1}}^{r}}+J_{4}(t),  \label{c-c(2)}
\end{eqnarray}%
\vspace{-0.6cm}
\begin{eqnarray}
t^{\mu _{r}}\,\Vert \nabla (v-\tilde{v})(t)\Vert _{\mathcal{M}_{r_{1}}^{r}}
&\leq &t^{\mu _{r}}\,\Vert \nabla e^{-\gamma t}e^{t\Delta }(v_{0}-\tilde{v}%
_{0})\Vert _{\mathcal{M}_{r_{1}}^{r}}  \notag \\
&&+t^{\mu _{r}}\int_{0}^{t}\Vert \nabla e^{-\gamma (t-\tau )}e^{(t-\tau
)\Delta }(u\cdot \nabla v+n-\tilde{u}\cdot \nabla \tilde{v}-\tilde{n})(\tau
)\Vert _{\mathcal{M}_{r_{1}}^{r}}\,d\tau  \notag \\
&\coloneqq&t^{\mu _{r}}\,\Vert \nabla e^{-\gamma t}e^{t\Delta }(v_{0}-\tilde{%
v}_{0})\Vert _{\mathcal{M}_{r_{1}}^{r}}+J_{5}(t)  \label{v-v}
\end{eqnarray}%
\vspace{-0.2cm} and \vspace{-0.2cm}
\begin{eqnarray}
t^{\mu _{p}}\,\Vert (u-\tilde{u})(t)\Vert _{\mathcal{M}_{p_{1}}^{p}} &\leq
&t^{\mu _{p}}\,\Vert e^{t\Delta }(u_{0}-\tilde{u}_{0})\Vert _{\mathcal{M}%
_{p_{1}}^{p}}+t^{\mu _{p}}\displaystyle\int_{0}^{t}\Vert e^{(t-\tau )\Delta }%
\mathbb{P}(u\cdot \nabla u-\tilde{u}\cdot \nabla \tilde{u})(\tau )\Vert _{%
\mathcal{M}_{p_{1}}^{p}}\,d\tau  \notag \\
&&+\,t^{\mu _{p}}\displaystyle\int_{0}^{t}\Vert e^{(t-\tau )\Delta }\mathbb{P%
}(nf-\tilde{n}f)(\tau )\Vert _{\mathcal{M}_{p_{1}}^{p}}\,d\tau  \notag \\
&\coloneqq&t^{\mu _{p}}\,\Vert e^{t\Delta }(u_{0}-\tilde{u}_{0})\Vert _{%
\mathcal{M}_{p_{1}}^{p}}+J_{6}(t)+J_{7}(t).  \label{u-u}
\end{eqnarray}

In view of \eqref{estB^2_4,2}, \eqref{estB^2_1,2}, \eqref{estDeltaB^2_4,2}, %
\eqref{estDeltaB^2_1,2}, \eqref{estB^3_4,3}, \eqref{estB^4_4,4}, %
\eqref{estL_3} and \eqref{estL_4}, we have the following estimates for the
integrals $J_{3}$, $J_{4}$, $J_{5}$, $J_{6}$ and $J_{7}$:
\begin{eqnarray}
J_{3}(t) &\leq &\tilde{C_{4,1}}\,\int_{0}^{1}(1-z)^{\mu _{p}-1}\,z^{-\mu
_{p}}\big((tz)^{\mu _{p}}\,\Vert (u-\tilde{u})(tz)\Vert _{\mathcal{M}%
_{p_{1}}^{p}}\,\Vert c\Vert _{X_{2}}  \notag \\
&&\text{ \ \ \ \ \ \ }+\Vert \tilde{u}\Vert _{X_{4}}\,\Vert (c-\tilde{c}%
)(tz)\Vert _{L_{\infty }}\big)\,dz  \notag \\
&&+\,\tilde{C_{5,1}}\,\int_{0}^{1}(1-z)^{l_{q}-1}\,z^{-l_{q}}\big(%
(tz)^{l_{q}}\,\Vert (n-\tilde{n})(tz)\Vert _{\mathcal{M}_{q_{1}}^{q}}\,\Vert
c\Vert _{X_{2}}  \notag \\
&&\text{ \ \ \ \ \ \ \ \ \ \ }+\Vert \tilde{n}\Vert _{X_{1}}\,\Vert (c-%
\tilde{c})(tz)\Vert _{L_{\infty }}\big)\,dz,\,\,\,\,\,\,\,\,\,  \label{J3}
\end{eqnarray}%
\vspace{-0.6cm}
\begin{eqnarray}
J_{4}(t) &\leq &\tilde{C_{4,2}}\,\int_{0}^{1}(1-z)^{\mu _{p}-1}\,z^{-\mu
_{p}-\mu _{r}}\,\big((tz)^{\mu _{p}}\,\Vert (u-\tilde{u})(tz)\Vert _{%
\mathcal{M}_{p_{1}}^{p}}\,\Vert c\Vert _{X_{2}}\,  \notag \\
&&\text{ \ \ \ \ \ \ }+\,\Vert \tilde{u}\Vert _{X_{4}}\,(tz)^{\mu
_{r}}\,\Vert \nabla (c-\tilde{c})(tz)\Vert _{\mathcal{M}_{r_{1}}^{r}}\,\big)%
dz  \notag \\
&&+\,\tilde{C_{5,2}}\,\int_{0}^{1}(1-z)^{l_{q}-\mu _{r}-1}\,z^{-l_{q}}\,\big(%
(tz)^{l_{q}}\Vert \,(n-\tilde{n})(tz)\Vert _{\mathcal{M}_{q_{1}}^{q}}\Vert
c\Vert _{X_{2}}  \notag \\
&&\text{ \ \ \ \ \ \ \ \ \ \ \ }+\Vert \tilde{n}\Vert _{X_{1}}\,\Vert (c-%
\tilde{c})(tz)\Vert _{L_{\infty }}\,\big)dz,  \label{J4}
\end{eqnarray}%
\vspace{-0.6cm}
\begin{eqnarray}
J_{5}(t) &\leq &\tilde{C_{6}}\,\int_{0}^{1}(1-z)^{\mu _{p}-1}\,z^{-\mu
_{p}-\mu _{r}}\big((tz)^{\mu _{p}}\,\Vert (u-\tilde{u})(tz)\Vert _{\mathcal{M%
}_{p_{1}}^{p}}\,\Vert v\Vert _{X_{3}}  \notag \\
&&\text{ \ \ \ \ \ \ }+\Vert \tilde{u}\Vert _{X_{4}}\,(tz)^{\mu _{r}}\Vert
\nabla (v-\tilde{v})(tz)\Vert _{\mathcal{M}_{r_{1}}^{r}}\,\big)dz  \notag \\
&&+\,\tilde{\alpha}\,\int_{0}^{1}(1-z)^{l_{q}-\mu
_{r}-1}\,z^{-l_{q}}\,(tz)^{l_{q}}\Vert (n-\tilde{n})(tz)\Vert _{\mathcal{M}%
_{q_{1}}^{q}}\,dz,  \label{J5}
\end{eqnarray}%
\vspace{-0.6cm}
\begin{eqnarray}
J_{6}(t) &\leq &\tilde{C_{7}}\,\int_{0}^{1}(1-z)^{\mu _{p}-1}\,z^{-2\mu
_{p}}\,\big((tz)^{\mu _{p}}\,\Vert (u-\tilde{u})(tz)\Vert _{\mathcal{M}%
_{p_{1}}^{p}}\,\Vert u\Vert _{X_{4}}  \notag \\
&&+\Vert \tilde{u}\Vert _{X_{4}}\,(tz)^{\mu _{p}}\,\Vert (u-\tilde{u}%
)(tz)\Vert _{\mathcal{M}_{p_{1}}^{p}}\,\big)dz\,\,\,\,\,\,\,\,\,\,\,\,\,
\label{J6}
\end{eqnarray}%
\vspace{-0.2cm} and \vspace{-0.2cm}
\begin{equation}
J_{7}(t)\leq \tilde{\beta}\,\int_{0}^{1}(1-z)^{l_{q}-\mu
_{p}-1}\,z^{-l_{q}}(tz)^{l_{q}}\Vert (n-\tilde{n})(tz)\Vert _{\mathcal{M}%
_{q_{1}}^{q}}\,dz.  \label{J7}
\end{equation}

Now we define
\begin{eqnarray*}
A_{1} &\coloneqq&\limsup_{t\rightarrow \infty }\,t^{l_{q}}\,\Vert n(\cdot
,t)-\tilde{n}(\cdot ,t)\Vert _{\mathcal{M}_{q_{1}}^{q}}, \\
&& \\
A_{2} &\coloneqq&\limsup_{t\rightarrow \infty }\,\Vert c(\cdot ,t)-\tilde{c}%
(\cdot ,t)\Vert _{L^{\infty }}, \\
&& \\
A_{3} &\coloneqq&\limsup_{t\rightarrow \infty }\,t^{\mu _{r}}\,\Vert \nabla
(c(\cdot ,t)-\tilde{c}(\cdot ,t))\Vert _{\mathcal{M}_{r_{1}}^{r}}, \\
&& \\
A_{4} &\coloneqq&\limsup_{t\rightarrow \infty }\,t^{\mu _{r}}\,\Vert \nabla
(v(\cdot ,t)-\tilde{v}(\cdot ,t))\Vert _{\mathcal{M}_{r_{1}}^{r}}, \\
&& \\
A_{5} &\coloneqq&\limsup_{t\rightarrow \infty }\,t^{\mu _{p}}\,\Vert u(\cdot
,t)-\tilde{u}(\cdot ,t)\Vert _{\mathcal{M}_{p_{1}}^{p}}.
\end{eqnarray*}%
Since $\left\Vert (n,c,v,u)\right\Vert _{\mathcal{X}},\left\Vert (\tilde{n},%
\tilde{c},\tilde{v},\tilde{u})\right\Vert _{\mathcal{X}}\leq
2K_{1}\varepsilon $, we have that $A_{1}$, $A_{2},\,A_{3},\,A_{4},\,A_{5}<%
\infty $. Taking the $\displaystyle\limsup_{t\rightarrow \infty }$ in %
\eqref{n-n}, \eqref{c-c(1)}, \eqref{c-c(2)}, \eqref{v-v} and \eqref{u-u},
and using \eqref{J1}, \eqref{J2} and \eqref{J3}-\eqref{J7}, we obtain

\begin{eqnarray}
A_{1} &\leq &0+\tilde{C_{1}}\,2\,K_{1}\,\varepsilon
\,\int_{0}^{1}\,(1-z)^{\mu _{p}-1}\,z^{-\mu _{p}-l_{q}}\,dz\,\left(
A_{5}+A_{1}\right)  \notag \\
&&+\tilde{C_{2}}\,2\,K_{1}\,\varepsilon \,\int_{0}^{1}\,(1-z)^{\mu
_{r}-1}\,z^{-l_{q}-\mu _{r}}\,\left( A_{1}+A_{3}\right) dz  \notag \\
&&+\,\tilde{C_{3}}\,2\,K_{1}\,\varepsilon \,\int_{0}^{1}\,(1-z)^{\mu
_{r}-1}\,z^{-l_{q}-\mu _{r}}\,\left( A_{1}+A_{4}\right) dz  \notag \\
&\leq &\,2\,K_{1}\,\varepsilon \left[ C_{1}\,\left( A_{1}+A_{5}\right)
+C_{2}\,\left( A_{1}+A_{3}\right) +C_{3}\,\left( A_{1}+A_{4}\right) \right] ,
\label{A1-aux} \\
A_{2} &\leq &0+\tilde{C_{4,1}}\,2\,K_{1}\,\varepsilon
\,\int_{0}^{1}\,(1-z)^{\mu _{p}-1}\,z^{-\mu _{p}}\,dz\,\,\left(
A_{5}+A_{2}\right)  \notag \\
&&+\tilde{C_{5,1}}\,2\,K_{1}\,\varepsilon
\,\int_{0}^{1}\,(1-z)^{l_{q}-1}\,z^{-l_{q}}\,dz\,\,\left( A_{1}+A_{2}\right)
\notag \\
&\leq &\,2\,K_{1}\,\varepsilon \left[ C_{4,1}\left( A_{2}+A_{5}\right)
+C_{5,1}\left( A_{1}+A_{2}\right) \right] ,  \notag \\
A_{3} &\leq &0+\tilde{C_{4,2}}\,2\,K_{1}\,\varepsilon
\,\int_{0}^{1}\,(1-z)^{\mu _{p}-1}\,z^{-\mu _{p}-\mu _{r}}\,dz\,\,\left(
A_{5}+A_{3}\right)  \notag \\
&&+\tilde{C_{5,2}}\,2\,K_{1}\,\varepsilon \,\int_{0}^{1}\,(1-z)^{l_{q}-\mu
_{r}-1}\,z^{-l_{q}}\,dz\,\,\left( A_{1}+A_{2}\right)  \notag \\
&\leq &2\,K_{1}\,\varepsilon \left[ C_{4,2}\,\left( A_{3}+A_{5}\right)
+C_{5,2}\,\left( A_{1}+A_{2}\right) \right] ,  \notag \\
A_{4} &\leq &0+\tilde{C_{6}}\,2\,K_{1}\,\varepsilon
\,\int_{0}^{1}\,(1-z)^{\mu _{p}-1}\,z^{-\mu _{p}-\mu _{r}}\,dz\,\left(
A_{5}+A_{4}\right)  \notag \\
&&+\tilde{\alpha}\,\int_{0}^{1}\,(1-z)^{l_{q}-\mu
_{r}-1}\,z^{-l_{q}}\,dz\,A_{1}  \notag \\
&\leq &2\,K_{1}\,\varepsilon \left[ C_{6}\,\left( A_{4}+A_{5}\right) \right]
+\alpha \,A_{1},  \notag \\
A_{5} &\leq &0+\tilde{C_{7}}\,2\,K_{1}\,\varepsilon
\,\int_{0}^{1}\,(1-z)^{\mu _{p}-1}\,z^{-2\mu _{p}}\,dz\,\left(
A_{5}+A_{5}\right)  \notag \\
&&+\tilde{\beta}\,\int_{0}^{1}\,(1-z)^{l_{q}-\mu
_{p}-1}\,z^{-l_{q}}\,dz\,A_{1}  \notag \\
&\leq &2\,K_{1}\,\varepsilon \left[ C_{7}\,2\,A_{5}\right] +\beta \,A_{1},
\notag
\end{eqnarray}%
where $%
\{C_{1},C_{2},C_{3},C_{4}=C_{4,1}+C_{4,2},C_{5}=C_{5,1}+C_{5,2},C_{6},C_{7}%
\} $ and $\{\alpha ,\beta \}$ are as in Lemmas \ref{LemaEstimativaBilinear}
and \ref{LemaEstimativaLinear}, respectively.

Recalling that $K_{1}=1+\alpha +\beta $ and $K_{2}=(\alpha +\beta
)(C_{1}+C_{2}+C_{3})+\sum_{i=1}^{7}C_{i}$ (see (\ref{K1eK2})) and summing
all $A_{i}$'s, we arrive at

\begin{eqnarray*}
A_{1}+A_{2}+A_{3}+A_{4}+A_{5} &\leq &2\,K_{1}\,\varepsilon \,\bigg[%
A_{1}\,\left( C_{1}+C_{2}+C_{3}+C_{5,1}+C_{5,2}\right) +A_{2}\,\left(
C_{4,1}+C_{5,1}+C_{5,2}\right) \\
&&+\,A_{3}\,\left( C_{2}+C_{4,2}\right) +A_{4}\,\left( C_{3}+C_{6}\right) \\
&&+A_{5}\,\left( C_{1}+C_{4,1}+C_{4,2}+C_{6}+2C_{7}\right) \bigg]+(\alpha
+\beta )A_{1} \\
&\leq &2\,K_{1}\,\varepsilon \,\bigg[A_{1}\,\left(
C_{1}+C_{2}+C_{3}+C_{5,1}+C_{5,2}+(\alpha +\beta )[C_{1}+C_{2}+C_{3}]\right)
\\
&&+A_{2}\,\left( C_{4,1}+C_{5,1}+C_{5,2}\right) +\,A_{3}\,\left(
C_{2}+C_{4,2}+(\alpha +\beta )C_{2}\right) \\
&&+A_{4}\,\left( C_{3}+C_{6}+(\alpha +\beta )C_{3}\right) \\
&&+A_{5}\,\left( C_{1}+C_{4,1}+C_{4,2}+C_{6}+2C_{7}+(\alpha +\beta
)C_{1}\right) \bigg]\text{ (by }\eqref{A1-aux}\text{)} \\
&\leq &2\,K_{1}\,\varepsilon \left( A_{1}+A_{2}+A_{3}+A_{4}+A_{5}\right)
\times \\
&&\bigg[C_{1}+C_{2}+C_{3}+C_{4,1}+C_{5,1}+C_{4,2}+C_{5,2}+C_{6} \\
&&+2C_{7}+(\alpha +\beta )(C_{1}+C_{2}+C_{3})\bigg].
\end{eqnarray*}%
As $C_{4}=C_{4,1}+C_{4,2}$ and $C_{5}=C_{5,1}+C_{5,2},$ note that $%
C_{1}+C_{2}+C_{3}+C_{4,1}+C_{5,1}+C_{4,2}+C_{5,2}+C_{6}+2C_{7}+(\alpha
+\beta )(C_{1}+C_{2}+C_{3})\leq 2\,K_{2}$, and then
\begin{equation*}
A_{1}+A_{2}+A_{3}+A_{4}+A_{5}\leq 4\,K_{1}\,K_{2}\,\varepsilon \,\left(
A_{1}+A_{2}+A_{3}+A_{4}+A_{5}\right) .
\end{equation*}%
Since $4K_{1}K_{4}\varepsilon <1$, it follows that $%
A_{1}=A_{2}=A_{3}=A_{4}=A_{5}=0$.

Now we turn to show that \eqref{volta} implies \eqref{ida}. We proceed as in
the estimates \eqref{n-n} and \eqref{c-c(1)}-\eqref{u-u} and use the
hypothesis $A_{1}=A_{2}=A_{3}=A_{4}=A_{5}=0$ (see (\ref{volta})) in order to
obtain
\begin{eqnarray*}
\lim_{t\rightarrow \infty }\,\sup \,t^{l_{q}}\Vert e^{t\Delta }(n_{0}-\tilde{%
n}_{0})\Vert _{\mathcal{M}_{q_{1}}^{q}} &\leq &A_{1}+\lim_{t\rightarrow
\infty }\sup \left( J_{1}(t)+J_{2}(t)\right) \\
&\leq &A_{1}+2\,K_{1}\,\varepsilon \,C_{1}\,\left( A_{1}+A_{5}\right) \\
&&+2\,K_{1}\,\varepsilon \,C_{2}\,\left( A_{1}+A_{3}\right)
+2\,K_{1}\,\varepsilon \,C_{3}\,\left( A_{1}+A_{4}\right) \\
&=&0+0+0+0=0, \\
\lim_{t\rightarrow \infty }\,\sup \,\Vert e^{t\Delta }(c_{0}-\tilde{c}%
_{0})\Vert _{L^{\infty }} &\leq &A_{2}+\lim_{t\rightarrow \infty }\sup
J_{3}(t) \\
&\leq &A_{2}+2\,K_{1}\,\varepsilon \,C_{4,1}\left( A_{2}+A_{5}\right)
+2\,K_{1}\,\varepsilon \,C_{5,1}\left( A_{1}+A_{2}\right) \\
&=&0+0+0=0, \\
\lim_{t\rightarrow \infty }\,\sup \,t^{\mu _{r}}\Vert \nabla e^{t\Delta
}(c_{0}-\tilde{c}_{0})\Vert _{\mathcal{M}_{r_{1}}^{r}} &\leq
&A_{3}+\lim_{t\rightarrow \infty }\sup J_{4}(t) \\
&\leq &A_{3}+2\,K_{1}\,\varepsilon \,C_{4,2}\,\left( A_{3}+A_{5}\right)
+2\,K_{1}\,\varepsilon \,C_{5,2}\,\left( A_{1}+A_{2}\right) \\
&=&0+0+0=0, \\
\lim_{t\rightarrow \infty }\,\sup \,t^{\mu _{r}}\Vert \nabla e^{-\gamma
t}e^{t\Delta }(v_{0}-\tilde{v}_{0})\Vert _{\mathcal{M}_{r_{1}}^{r}} &\leq
&A_{4}+\lim_{t\rightarrow \infty }\sup J_{5}(t) \\
&\leq &A_{4}+2\,K_{1}\,\varepsilon \,C_{6}\,\left( A_{4}+A_{5}\right)
+\alpha \,A_{1} \\
&=&0+0+0=0
\end{eqnarray*}%
and
\begin{eqnarray*}
\lim_{t\rightarrow \infty }\,\sup \,t^{\mu _{p}}\Vert e^{t\Delta }(u_{0}-%
\tilde{u}_{0})\Vert _{\mathcal{M}_{p_{1}}^{p}} &\leq
&A_{5}+\lim_{t\rightarrow \infty }\sup \left( J_{6}(t)+J_{7}(t)\right) \\
&\leq &A_{5}+2\,K_{1}\,\varepsilon \,C_{7}\,2A_{5}+\beta \,A_{1} \\
&=&0+0+0=0,
\end{eqnarray*}%
and we are done.\begin{flushright}$\blacksquare$\end{flushright}

\noindent\textbf{{\large {Acknowledgments.}}} LCF Ferreira was supported by
FAPESP and CNPq, Brazil. M. Postigo was supported by CAPES and CNPq, Brazil.

\end{document}